%%%%%%%%%%%
% typoref.tex. V : January 18, 2000.
% Author : Anthony PHAN
% Warning : syntaxe +- LaTeX
% Sources :
% T. Lachand--Robert, ``La Ma\^\i trise de \TeX'',
% R\'ef\'erences crois\'ees;
% latex.ltx's sources;
% and of course the \TeX book.
%%%%%%%%%%%%%%%%%%%%%%%%%%%%%%%%%%%%%%%%%%%%%%%%%%%%%
%

\catcode`\@=11

\magnification=1200
\baselineskip=14pt

\pretolerance=500    \tolerance=1000 \brokenpenalty=5000

\catcode`\;=\active
\def;{\relax\ifhmode\ifdim\lastskip>\z@
\unskip\fi\kern.2em\fi\string;}

\overfullrule=0mm

\catcode`\!=\active
\def!{\relax\ifhmode\ifdim\lastskip>\z@
\unskip\fi\kern.2em\fi\string!}

\catcode`\?=\active
\def?{\relax\ifhmode\ifdim\lastskip>\z@
\unskip\fi\kern.2em\fi\string?}

\frenchspacing

\newif\ifpagetitre            \pagetitretrue
\newtoks\hautpagetitre        \hautpagetitre={ }
\newtoks\baspagetitre         \baspagetitre={1}

\newtoks\auteurcourant        \auteurcourant={ }
\newtoks\titrecourant
\titrecourant={ }
\newtoks\hautpagegauche       \newtoks\hautpagedroite
\hautpagegauche={\hfill\sevenrm\the\auteurcourant\hfill}
\hautpagedroite={\hfill\sevenrm\the\titrecourant\hfill}

\newtoks\baspagegauche       \baspagegauche={\hfill\rm\folio\hfill}

\newtoks\baspagedroite       \baspagedroite={\hfill\rm\folio\hfill}

\headline={
\ifpagetitre\the\hautpagetitre
\global\pagetitrefalse
\else\ifodd\pageno\the\hautpagedroite
\else\the\hautpagegauche\fi\fi}

\footline={\ifpagetitre\the\baspagetitre
\global\pagetitrefalse
\else\ifodd\pageno\the\baspagedroite
\else\the\baspagegauche\fi\fi}

\def\date{\ {\the\day}\
\ifcase\month\or Janvier\or F\'evrier\or Mars\or Avril
\or Mai \or Juin\or Juillet\or Ao\^ut\or Septembre
\or Octobre\or Novembre\or D\'ecembre\fi\
{\the\year}}

\def\up#1{\raise 1ex\hbox{\sevenrm#1}}

\def\cqfd{\unskip\kern 6pt\penalty 500
\raise -2pt\hbox{\vrule\vbox to 10pt{\hrule width 4pt
\vfill\hrule}\vrule}\par\medskip}

\def\section#1{\vskip 7mm plus 20mm minus 1.5mm\penalty-50
\vskip 0mm plus -20mm minus 1.5mm\penalty-50
{\bf\noindent#1}\nobreak\smallskip}

\def\subsection#1{\medskip{\bf#1}\nobreak\smallskip}

\def\displaylinesno #1{\dspl@y\halign{
\hbox to\displaywidth{$\@lign\hfil\displaystyle##\hfil$}&
\llap{$##$}\crcr#1\crcr}}

\def\ldisplaylinesno #1{\dspl@y\halign{
\hbox to\displaywidth{$\@lign\hfil\displaystyle##\hfil$}&
\kern-\displaywidth\rlap{$##$}
\tabskip\displaywidth\crcr#1\crcr}}

\def\hfl#1#2{\smash{\mathop{\hbox to 12 mm{\rightarrowfill}}
\limits^{\scriptstyle#1}_{\scriptstyle#2}}}

%
% style (look at the behavior of \item dans \bibitem too,
% and at one ,\  in \re@dreferenceslist)
% Feel free to change: 	\bibn@me (title like ``R\'ef\'erences'')
%			\bibliographym@rk (general style)
%
\def\bibn@me{R\'ef\'erences}
\def\bibliographym@rk{\centerline{{\sc\bibn@me}}
	\sectionmark\section{\ignorespaces}{\unskip\bibn@me}
	\bigbreak\bgroup
	\ifx\ninepoint\undefined\relax\else\ninepoint\fi}
%
% Beware of the \bgroup: it will be closed by \endthebibliography
%
% \refsp@ce is the spacing command that appens between multiple
% references.
%
\let\refsp@ce=\
\let\bibleftm@rk=[
\let\bibrightm@rk=]
%
% if you want more space between brackets...
%\let\refsp@ce=\thinspace
%\def\bibleftm@rk{[\thinspace}
%\def\bibrightm@rk{\thinspace]}
%
% frenchy stuff
%
\def\numero{n\raise.82ex\hbox{$\fam0\scriptscriptstyle
o$}~\ignorespaces}
%
% new variables
%
\newcount\equationc@unt
\newcount\bibc@unt
\newif\ifref@changes\ref@changesfalse
\newif\ifpageref@changes\ref@changesfalse
\newif\ifbib@changes\bib@changesfalse
\newif\ifref@undefined\ref@undefinedfalse
\newif\ifpageref@undefined\ref@undefinedfalse
\newif\ifbib@undefined\bib@undefinedfalse
\newwrite\@auxout
%
% mark an equation
%
%\def\eqnum{\global\advance\equationc@unt by 1%
%\edef\lastref{\number\equationc@unt}%
%\eqno{(\lastref)}}
%
% One can reference anything, just copy the former macro
% and use it so: \machin \label{truc}
% In machin you would have defined \lastref by some number
% or any text.
%
% References macros
%
% The next macros are the core of \ref and \cite commands.
% Its first argument may be ref, pageref or bib.
%
% It is too tricky to be explained.
% (It is a bit recursive.)
% It allows using \cite or \ref or ...
% with arbitrary many arguments,
% for instance:
% \cite{knuth1,knuth2,ma pomme}
%
% First argument is always ref, pageref or bib.
%
\def\re@dreferences#1#2{{%
	\re@dreferenceslist{#1}#2,\undefined\@@}}
\def\re@dreferenceslist#1#2,#3\@@{\def\next{#2}%
	\expandafter\ifx\csname#1@@\meaning\next\endcsname\relax
	??\immediate\write16
	{Warning, #1-reference "\next" on page \the\pageno\space
	is undefined.}%
	\global\csname#1@undefinedtrue\endcsname
	\else\csname#1@@\meaning\next\endcsname\fi
	\ifx#3\undefined\relax
	\else,\refsp@ce\re@dreferenceslist{#1}#3\@@\fi}
%
% notice that the former ``,\refsp@ce'' will separate
% multiple arguments. But beware of spaces
% while defining a reference or calling for it!
%
% tricky thing: \newlabel has two arguments
% {labelname}{{\lastref}{\pageref}}
% The second argument is read as two arguments
% by \newl@bel. This was necessary to get
% a jobname.aux containing the same syntax
% LaTeX would produce and use.
%
\def\newlabel#1#2{{\def\next{#1}\newl@bel#2}}
\def\newl@bel#1#2{%
	\expandafter\xdef\csname ref@@\meaning\next\endcsname{#1}%
	\expandafter\xdef\csname pageref@@\meaning\next\endcsname{#2}}
\def\label#1{{%
	\toks0={#1}\message{ref(\lastref) \the\toks0,}%
	\ignorespaces\immediate\write\@auxout%
	{\noexpand\newlabel{\the\toks0}{{\lastref}{\the\pageno}}}%
	\def\next{#1}%
	\expandafter\ifx\csname ref@@\meaning\next\endcsname\lastref%
	\else\global\ref@changestrue\fi%
	\newlabel{#1}{{\lastref}{\the\pageno}}}}
\def\ref#1{\re@dreferences{ref}{#1}}
\def\pageref#1{\re@dreferences{pageref}{#1}}
%
% bibliography macros
%
\def\bibcite#1#2{{\def\next{#1}%
	\expandafter\xdef\csname bib@@\meaning\next\endcsname{#2}}}
\def\cite#1{\bibleftm@rk\re@dreferences{bib}{#1}\bibrightm@rk}
%
% The argument of \beginthebibliography
% is any sequence of numerals which will represent
% the maximum \item's length. If you have less than 9
% \bibitem's, this argument may be {any numeral}.
% if you have between 100 and 999 \bibitem's
% this argument may be {any three numerals},
% and so on.
%
\def\beginthebibliography#1{\bibliographym@rk
	\setbox0\hbox{\bibleftm@rk#1\bibrightm@rk\enspace}
	\parindent=\wd0
	\global\bibc@unt=0
	\def\bibitem##1{\global\advance\bibc@unt by 1
		\edef\lastref{\number\bibc@unt}
		{\toks0={##1}
		\message{bib[\lastref] \the\toks0,}%
		\immediate\write\@auxout
		{\noexpand\bibcite{\the\toks0}{\lastref}}}
		\def\next{##1}%
		\expandafter\ifx
		\csname bib@@\meaning\next\endcsname\lastref
		\else\global\bib@changestrue\fi%
		\bibcite{##1}{\lastref}
		\medbreak
		\item{\hfill\bibleftm@rk\lastref\bibrightm@rk}%
		}
	}
\def\endthebibliography{\egroup\par}
%
% THE NEXT MACRO MUST BE INCLUDED
% IN THE \BYE COMMAND. FOR INSTANCE:
%
    %\catcode`@=11
    \outer\def\bye{\@closeaux
    	\par\vfill\supereject\end}
    %\catcode`@=12
%
\def\@closeaux{\closeout\@auxout
	\ifref@changes\immediate\write16
	{Warning, changes in references.}\fi
	\ifpageref@changes\immediate\write16
	{Warning, changes in page references.}\fi
	\ifbib@changes\immediate\write16
	{Warning, changes in bibliography.}\fi
	\ifref@undefined\immediate\write16
	{Warning, references undefined.}\fi
	\ifpageref@undefined\immediate\write16
	{Warning, page references undefined.}\fi
	\ifbib@undefined\immediate\write16
	{Warning, citations undefined.}\fi}
%
% initialization of jobname.aux
%
\immediate\openin\@auxout=\jobname.aux
\ifeof\@auxout \immediate\write16
     {Creating file \jobname.aux}
\immediate\closein\@auxout
\immediate\openout\@auxout=\jobname.aux
\immediate\write\@auxout {\relax}%
\immediate\closeout\@auxout
\else\immediate\closein\@auxout\fi
%
% Let's read this file and open it out
%
\input\jobname.aux \par
\immediate\openout\@auxout=\jobname.aux
% this file will be closed by {\bf y}e.
%
% That's all, folks!
%

    \def\Q{{\bf {Q}}}

\def\Z{{\bf Z}} 
\def\R{{\bf R}} 
\def\C{{\bf C}}

\def\bP{{\bf P}}

\def\bR{{\bf R}}
\def\bZ{{\bf Z}}

\def\bx{{\bf x}}
\def\by{{\bf y}}

  \def\pro{\noindent {\bf{Proof :
}}}
   \def\resp{{\rm resp.  }}

   \def\al{{\alpha}}
  \def\eps{{\varepsilon}}

\def\om{{\omega}}

\def\hom{{\hat{\omega}}}
\def\homega{{\hat{\omega}}}
\def\and{\quad\hbox{and}\quad}

 \def\ux{{\underline{x}}} \def\uy{{\underline
y}}

\def\build#1_#2^#3{\mathrel{\mathop{\kern 0pt#1}\limits_{#2}^{#3}}}

%\newfam\gothfam \scriptscriptfont\gothfam=\fivegoth
%\textfont\gothfam=\tengoth \scriptfont\gothfam=\sevengoth
%\def\goth{\fam\gothfam\tengoth}

%\def\cqfd{\unskip\kern 6pt\penalty 500 \raise 0pt\hbox{\vrule\vbox
%to6pt{\hrule width 6pt \vfill\hrule}\vrule}\par}

\def\pro{\noindent {\it Proof. }}
\def\proof{\noindent {\it Proof. }}

\def\smallsquare{\vbox{\hrule\hbox{\vrule height 1 ex\kern 1
ex\vrule}\hrule}}
\def\cqfd{\hfill \smallsquare\vskip 3mm}

\def\hw{{\hat w}}
\def\hla{{\hat \lambda}}
\def\hlambda{{\hat \lambda}}
\def\utheta{{\underline {\xi}}}
\def\uxi{{\underline{\xi}}}

\def\om{{\omega}}

\def\uu{{\uxi}}
\def\cC{{\cal C}}
\def\psibot{{\underline{\psi}}}
\def\psitop{{\overline{\psi}}}
\def\uL{{\bf L}}

\def\bibn@me{R\'ef\'erences bibliographiques}
%\input typpo
%
%\catcode`@=11
\def\bibliographym@rk{\bgroup}
%
% {\bf y}e est modifie pour la biblio et la table des matieres
%
\outer\def\bye{ 	\par\vfill\supereject\end}

%%%%%%%%%%%%%%%%%%%%%%%%%%%%%%%%%%%%%%%%%%
%\catcode`@=12

\null

%\centerline{}

\vskip 2mm

\hfill{\it Preliminary draft}

\vskip 6mm

\centerline{\bf Exponents of Diophantine approximation}

\vskip 8mm

\centerline{\it Yann Bugeaud \footnote{}{\rm
2000 {\it Mathematics Subject Classification : } 11J13.  
To be published in: `Dynamics and Number Theory', 
Proceedings of the Durham Easter School 2014.}
}

\vskip 11mm

We survey classical and recent results on exponents of Diophantine
approximation. We give only a few proofs and highlight several open 
problems. 

\vskip 5mm

\vskip 5mm

\centerline{\bf 1. Introduction and generalities}

\vskip 5mm

%This section is inspired by Section 1 of \cite{BuLa05b}.

Let $\xi$ be an irrational real number.
It follows from 
the theory of continued fractions that there are infinitely many rational
numbers $p/q$ with $q \ge 1$ and such that
$$
|q \xi - p| \le q^{-1}.   
$$
Said differently, for arbitrarily large integers $Q$, there exist
integers $p$ and $q$ with $1 \le q \le Q$ and $|q \xi - p| \le Q^{-1}$.
%$$|q \xi - p| \le Q^{-1}.  \eqno (1.1) $$
The Dirichlet Theorem asserts much more, namely that, for 
{\it every} integer $Q \ge 1$, there exist
integers $p$ and $q$ with $1 \le q \le Q$ and $|q \xi - p| \le Q^{-1}$.
% $$ |q \xi - p| \le Q^{-1}.  \eqno (1.2) $$

A question arises then naturally: is there some specific 
irrational real number $\xi$ for which it is possible to improve
the above statements, that is, to get the above inequalities with $Q^{-1}$
replaced by $Q^{-w}$ for some real number $w > 1$. 
This leads to the introduction of the exponents of approximation
$w_1$ and $\hw_1$. 

\proclaim Definition 1.1.
Let $\xi$ be a real number. 
We denote by $w_1(\xi)$ the supremum of the real numbers $w$
for which there exist arbitrarily large integers $Q$ 
and integers $p$ and $q$ with $1 \le q \le Q$ and
$$
|q \xi - p| \le Q^{-w}.  
$$
We denote by $\hw_1 (\xi)$ the supremum of the real numbers $\hw$
such that, for every sufficiently large integer $Q$, there are 
integers $p$ and $q$ with $1 \le q \le Q$ and
$$
|q \xi - p| \le Q^{-\hw}.   
$$

As observed by Khintchine \cite{Kh26b}, every irrational real number $\xi$
satisfies $\hw_1 (\xi) = 1$, since 
there are arbitrarily large integers $Q$ such that the inequality $|q \xi - p| \le 1/ (2Q)$
has no solutions in integers $p, q$ with $1 \le q \le Q$; see 
the proof of Proposition 2.4 below. 
However, for any $w>1$, there exist irrational real numbers $\xi$
such that, for {\it arbitrarily large} integers $Q$, the equation
$$
|q \xi - p| \le Q^{-w}
$$
has a solution in integers $p$ and $q$ with $1 \le q \le Q$. It suffices for instance
to take $\xi = \sum_{j \ge 1} \, 2^{-\lfloor (w+1)^j \rfloor}$, where $\lfloor \cdot \rfloor$ denotes
the integer part. 

Obviously, the quality of approximation strongly depends upon
whether we are interested in a uniform statement (that is, a statement valid
for every $Q$, or for every sufficiently large $Q$) or in a statement
valid only for arbitrarily large values of $Q$.

The general framework is the following.
For any (column) vector ${{\underline \theta}}$   in $\bR^n$, we denote by $\vert
{{\underline \theta}}\vert$
the  maximum of the   absolute values of its coordinates  and by
$$
\|  {{\underline \theta}}\|  = \min_{\uy \in \bZ^n} \vert {{\underline \theta}}  - \uy\vert
$$
the maximum of the distances of its  coordinates to the  rational
integers.

\proclaim Definition 1.2.
Let $n$ and $m$ be positive integers and $A$ a real matrix
with $n$ rows and $m$ columns.
Let ${{\underline \theta}}$ be an $n$-tuple 
of real numbers. We denote by  $w(A,{{\underline \theta}})$
the supremum of the real numbers $w$ for which, {\it for arbitrarily
large real numbers $X$}, the inequalities
$$
\|  A\ux - {{\underline \theta}}\|  \le X^{-w} \and \vert \ux \vert \le X
\eqno{(1.1)}
$$
have a solution  $\ux$ in $\bZ^m$. 
We denote by 
$\hat{w}(A,{{\underline \theta}})$ the supremum of the real numbers $w$
for which, {\it for all sufficiently large positive real numbers  $X$},
the inequalities (1.1) have an integer solution $\ux$ in $\bZ^m$.

In the sequel, we consistently use the symbol $\, \hat{ }$
to indicate that we require a uniform existence of solutions.

The lower bounds
$$
w(A,{{\underline \theta}}) \ge \hat{w}(A,{{\underline \theta}}) \ge 0
$$
are then obvious. We define furthermore two homogeneous  exponents
$w(A)$ and $\hat{w}(A)$ as in (1.1) with ${{\underline \theta}} = {}^t (0,\dots , 0)$,
requiring moreover  that the integer solution $\ux$  should be
non-zero.
The  uniform exponent $\hat{w}$
was first introduced and studied by Jarn\'\i k in the 1930s.

The transposed  matrix of a matrix
$A$ is denoted by ${}^t A$.
Furthermore, $1 / + \infty$ is understood to be $0$.
The following result, established in \cite{Cas,BuLa05b},
shows that the usual (resp. uniform) inhomogeneous exponents are
strongly related
to the uniform (resp. usual) homogeneous exponents.

\proclaim Theorem 1.3.
Let $n$ and $m$ be positive integers and $A$ a real matrix
with $n$ rows and $m$ columns.
For any $n$-tuple ${{\underline \theta}}$ of real numbers, we have the lower bounds
$$
w(A,{{\underline \theta}}) \ge {1\over \hat{w}({}^tA)} \and
\hat{w}(A,{{\underline \theta}}) \ge {1\over {w}({}^tA)} ,
\eqno(1.2)
$$
with equality  in (1.2) for almost all ${{\underline \theta}}$ with respect to the
Lebesgue measure on $\bR^n$.

If  the subgroup $G= {}^tA\bZ^n + \bZ^m$ of $\bR^m$ generated
by the $n$ rows of the matrix $A$ together with $\bZ^m$
has maximal rank $m+n$,  then Kronecker's Theorem
asserts that the dual subgroup $\Gamma = A \bZ^m +\bZ^n$ of $\bR^n$
generated by the
$m$ columns of  $A$ and by $\bZ^n$  is dense in
$\bR^n$. In this respect,  Theorem 1.3 may be viewed as a  measure of
the density of $\Gamma$.
In the case where the rank of  $G$ is
$< m+n$, we clearly have
$$
\hat{w}({}^t A) = w({}^t A) =+\infty \and
\hat{w}(A,{{\underline \theta}})={w}(A,{{\underline \theta}})=0, 
$$
for every $n$-tuple ${{\underline \theta}}$ located outside a discrete family of
parallel hyperplanes  in $\bR^n$. The
assertion of Theorem 1.3 is then obvious.

Cassels' book \cite{Cas} remains an invaluable reference for these and related questions.

In the sequel of the text, we restrict our attention 
to the cases where $A$ is either a row or a column
matrix: 
$$
A = (\xi_1, \dots , \xi_n)  \quad  \hbox{\rm or} \quad A =
{}^t(\xi_1, \dots , \xi_n).
$$
This amounts to consider small values of the linear form
$$
|x_0 +  x_1 \xi_1 + \ldots + x_n \xi_n|, \quad
\hbox{where $x_0, x_1, \ldots , x_n \in \Z$,}
$$
or simultaneous approximation to $\xi_1, \ldots , \xi_n$ by rational numbers
with the same denominator, that is, small values of the quantity
$$
\max_{1 \le j \le n} \, |x_0 \xi_j - x_j|, \quad
\hbox{where $x_0, x_1, \ldots , x_n \in \Z$.}
$$
Furthermore, among the elements $\uxi = (\xi_1, \dots , \xi_n)$ in $\R^n$, we 
mainly focus on the points
$$
(\xi, \xi^2, \ldots , \xi^n)
$$
whose coordinates are the $n$ first
successive powers of a real number $\xi$. However, some of the
results stated below hold for a general $n$-tuple $\uxi$, as will be indicated in due course. 
This is in particular the case in Sections 6 and 7. 

The present paper is organized as follows. 
In Section 2, we define six exponents of approximation attached to real numbers and
give their first properties. We discuss in Section 3 how these exponents are interrelated
and study the values taken by these exponents
in Sections 4 and 5. Intermediate exponents are introduced in Section 6. 
Parametric geometry of numbers, a deep and powerful, new theory introduced by 
Schmidt and Summerer \cite{ScSu09,ScSu13a} and developed by Roy \cite{Roy15a},
is briefly described in Section 7. 
Recent results on the existence of real numbers which are badly approximable by
algebraic numbers of bounded degree are discussed in Section 8.
The final section gathers several open problems and suggestion for further research. 

The notation $a \gg_d b$ means that $a$ exceeds $b$ times a constant 
depending only on $d$. When $\gg$ is written 
without any subscript, it means that the constant is absolute.
We write $a \asymp b$ if both
$a \gg b$ and $a \ll b$ hold. 

\vskip 5mm

\centerline{\bf 2. Further definitions and first results}

\vskip 5mm

Mahler \cite{Mah32} and Koksma \cite{Ko39} have introduced in the 1930s two classifications
of real numbers in terms of their properties of approximation
by algebraic numbers.
Here and below, the height $H(P)$ of a polynomial $P(X)$ is the maximum of the moduli of
its coefficients and the height $H(\alpha)$ of an algebraic number $\alpha$
is the height of its minimal polynomial over $\Z$.

\proclaim Definition 2.1.
Let $n$ be a positive integer.
Let $\xi$ be a real number. 
We denote by $w_n(\xi)$ the supremum of the real numbers $w$ for which
the inequality
$$
0 < |P(\xi)| \le H(P)^{-w}
$$
is  satisfied for infinitely many polynomials
$P(X)$ with integer coefficients and degree at most $n$.
We denote by $w_n^*(\xi)$ the supremum of the real numbers $w^*$ for which 
the inequality
$$
0 < |\xi - \alpha| \le H(\alpha)^{- w^* -1}  \eqno (2.1)
$$
is  satisfied for infinitely many algebraic numbers
$\alpha$ of degree at most $n$.

The $-1$ in the exponent of the right hand side of (2.1) can be
explained as follows. 
Let $P(X)$ be an integer polynomial of degree $n \ge 2$ and $\xi$
be a real number not algebraic of degree at most $n$. 
Observe that $P(\xi) P'(\xi) \not= 0$ and
$$
{P'(\xi) \over P(\xi)} = \sum_{\alpha : P(\alpha) = 0} \, {1 \over \xi - \alpha}.
$$
Consequently, if $\alpha$ is the root of $P(X)$ which is closest to $\xi$, then we have 
$$
|\xi - \alpha| \le n |P(\xi) / P'(\xi)|.
$$
Since $|P'(\xi)|$ is 
expected to be roughly equal to
$H(P)$ (this is the case unless $P(X)$ has two roots close to $\xi$), 
we expect the estimation 
$$
|\xi - \alpha| \ll_n |P(\xi)|  \cdot H(P)^{-1}. 
$$ 
The latter upper bound does not hold when $|P'(\xi)|$ is small, that is, when $P(X)$
has two or more roots close to $\xi$. 

The behaviour of the sequences $(w_n (\xi))_{n \ge 1}$ and
$(w_n^*(\xi))_{n \ge 1}$ determines the localisation of $\xi$ in
Mahler's and Koksma's classifications, respectively (see Chapter 3 of \cite{BuLiv});
however, the exact determination of $w_n (\xi)$ and $w_n^*(\xi)$
for a specific real number $\xi$ is usually extremely difficult.

%The Dirichlet box principle (or, equivalently,
%Minkowski's theorem) readily implies that $w_n(\xi) \ge n$ for
%any positive integer $n$ and any real number $\xi$ not algebraic of degree $\le n$.
%It is a longstanding problem to decide whether the same
%result remains true for the exponents $w_n^*$.
%In particular, Davenport and
%Schmidt \cite{DaSc69} have investigated this question in a dual way,
%via  simultaneous rational approximation to the successive powers of $\xi$. 

We introduced in \cite{BuLa05a} four further exponents of approximation.
They implicitly appeared
previously in articles of Jarn\'\i k, Davenport and Schmidt, among others.

\proclaim Definition 2.2.
Let $n$ be a positive integer.
Let $\xi$ be a real number. 
We denote by $\lambda_n(\xi)$  
the supremum of the real numbers $\lambda$ such that
the inequality
$$
\qquad  \max_{1 \le j \le n} \,
|x_0 \xi^j - x_j| \le  |x_0|^{-\lambda}
$$
has infinitely  many solutions in integers $x_0, \dots,x_n$
with $x_0\not= 0$.

The  three exponents $w_n$, $w^*_n$ and $\lambda_n$
have the common feature to be defined by the existence
of infinitely  many solutions for some set of Diophantine inequalities.
We attach to them three exponents defined
by a condition of uniform existence of solutions.

\proclaim
Definition 2.3.
Let $n$ be a positive integer and let $\xi$ be a real number.
We denote by $\hw_n(\xi)$
the supremum of the real numbers $\hw$ such that, for
any sufficiently large real number $X$, the inequalities
$$
0 < |x_n \xi^n + \ldots + x_1 \xi + x_0| \le X^{-\hw}, \qquad  \max_{0
\le j \le n} \, |x_j| \le X,
$$
have a solution in integers $x_0, \ldots, x_n$.
We denote by  $\hw^*_n(\xi)$  the
  supremum of the real numbers $\hw^*$ such that, for
any sufficiently large real number $X$,
there exists an algebraic real number
  $\alpha$ with degree at most $n$ satisfying
$$
0 < |\xi - \alpha| \le H(\alpha)^{-1} \, X^{-\hw^*} \quad \hbox{and} \quad H(\alpha) \le X.
$$
We denote by $ \hat{\lambda}_n(\xi)$ the supremum  of the real numbers
$\hlambda$ such that, for
any sufficiently large real number $X$, the inequalities
$$
0 < |x_0| \le X, \qquad \max_{1 \le j \le n} \,
|x_0 \xi^j - x_j| \le X^{-\hlambda},
$$
have a solution in integers $x_0, \ldots, x_n$.

This survey is mainly devoted to an overview of general
results on the six exponents of approximation 
$w_n,w_n^*,\lambda_n, \hat{w}_n,\hat{w}_n^*$ and
$\hat{\lambda}_n$, whose values are connected through various inequalities.
Notice that   the exponents  `hat'  are uniformly
bounded in term of   $n$.  The exact determination of the upper bounds
is an important problem towards the Wirsing Conjecture (see Section 4) 
or related questions, such as the approximation 
of transcendental real numbers by algebraic integers. 
We begin with some easy properties.

Our first result goes back to Khintchine \cite{Kh26b}. 

\proclaim Proposition 2.4.
For any irrational real number $\xi$, we have
$$
\hat{w}_1 (\xi) = \hat{w}^*_1 (\xi) = \hat{\lambda}_1 (\xi) =1.
$$
for every irrational real number $\xi$. 

\pro
Let $\xi$ be an irrational real number and $(p_\ell / q_\ell)_{\ell \ge 1}$ denote 
the sequence of its convergents. 
Let $\ell \ge 4$ and $q$ be integers with $1 \le q \le q_{\ell} - 1$.
Observe that $q_{\ell} - q_{\ell - 1} \ge q_{\ell - 2} \ge q_2 \ge 2$. 
Then, it follows from the theory of continued fractions that 
$$
\| q \xi \| \ge \|q_{\ell-1} \xi\| \ge {1 \over  q_\ell + q_{\ell - 1}} \ge {1 \over 2 (q_{\ell} - 1)}.
$$
This shows that, setting $Q = q_{\ell} - 1$, the inequality
$|q \xi - p| \le 1/ (2 Q)$ has no solutions in integers $p, q$ with $1\le q \le q_{\ell} - 1$.
Consequently, $\hw_1 (\xi)$ is equal to $1$. 
\cqfd

We gather in the next theorem several easy results on 
our six classical exponents of approximation.

\proclaim Theorem 2.5. 
For any positive integer $n$ and any real number
 $\xi$ which is not algebraic of degree
$\le n$, we have
$$
n \le \hat{w}_n (\xi) \le w_n (\xi), \qquad
{1\over n} \le\hat{\lambda}_n(\xi) \le \min\{1,\lambda_n(\xi)\}
$$
and
$$
1 \le \hat{w}^*_n (\xi) \le \min\{ w_n^* (\xi) ,\hat{w}_n (\xi)\}
\le \max\{ w_n^* (\xi) ,\hat{w}_n (\xi)\}\le w_n (\xi). \eqno (2.2)
$$

Corollary 5.4 shows that the lower bound $\hat{w}_n^*(\xi) \ge 1$ in (2.2) is 
best possible.

%if $\xi$ is a Liouville number, then $\hat{w}^*_n (\xi) = 1$
%holds for any positive integer $n$.

\proof 
Let $n$ and $\xi$ be as in the statement of the theorem. 
Let $\alpha$ be a real algebraic number close to $\xi$ and $P_{\alpha}(X)$ the
minimal defining polynomial of $\alpha$ over $\Z$. Then, we have
$|P_\alpha (\xi)| \ll_n H(P_{\alpha}) \cdot |\xi - \alpha|$ (see e.g. Lemma A.6 of \cite{BuLiv})
and we deduce the upper bound $w_n^* (\xi) \le w_n (\xi)$. 
The same argumentation gives also that $\hat{w}_n^*(\xi) \le \hat{w}_n (\xi)$.
The upper bounds $\hat{w}_n (\xi) \le w_n (\xi)$, $\hat{w}^*_n (\xi) \le w_n^* (\xi)$ 
and $\hat{\lambda}_n (\xi) \le \lambda_n (\xi)$
are consequences of the definitions, 
while the lower bounds $ \hat{w}_n (\xi)\ge n$ and
$\hat{\lambda}_n (\xi) \ge 1/n$
follow from Dirichlet's box principle (or from Minkowski's theorem).
Moreover, we obviously have
$\hat{w}^*_n (\xi) \ge \hat{w}^*_1 (\xi)=1$
and
$\hat{\lambda}_n (\xi) \le \hat{\lambda}_1 (\xi) =1$. \cqfd

The next theorem was pointed out in \cite{BuLa05a}. 

\proclaim Theorem 2.6.  
For any positive integer $n$ and any real number $\xi$ not algebraic of 
degree at most $n$, we have
$$
\hat{w}_n^* (\xi) \ge {w_n (\xi) \over w_n (\xi) -n +1}  \eqno (2.3)
$$
and
$$
w_n^* (\xi) \ge {\hat{w}_n (\xi) \over \hat{w}_n (\xi) -n +1}  \eqno (2.4)
$$

Wirsing \cite{Wir} proved a weaker version of (2.3) in which the left-hand side is replaced
by the quantity ${w}_n^* (\xi)$. His result is also weaker than (2.4), since 
$\hat{w}_n (\xi) \ge w_n(\xi)$. 

\proof 
We follow an argumentation of Wirsing \cite{Wir}.
Let $n \ge 2$ and $\xi$ be as in the statement of the theorem. 
We first establish (2.3).
If $w_n (\xi)$ is infinite, then (2.3) reduces to $\hw_n^* (\xi) \ge 1$, a statement
established in Theorem 2.5.
Assume that $w_n (\xi)$ is finite. 
Let $\eps > 0$ and set $w = w_n(\xi) (1 + \eps)^2$. 
Minkowski's theorem implies that there exist a
constant $c$ and, for any positive real number $H$, 
a non zero integer polynomial
$P(X)$ of degree at most $n$ such that
$$
|P(\xi)| \le H^{-w}, \quad |P(1)|, \ldots, |P(n-1)| \le H \quad
{\rm and} \quad |P(n)| \le c H^{w-n+1}. \eqno (2.5)
$$
The definition of  $w_n(\xi)$ and the first inequality of (2.5) show 
that $H(P) > H^{1 + \eps}$ if $H$ is large enough.
Consequently, $P(X)$ has some (necessarily
real) root in the neighbourhood of each of the points $\xi$, $1, \ldots, 
n-1$. Denoting by $\alpha$
its closest root to $\xi$ and recalling that
$H(\alpha) \ll_n H(P)$ (see for example Lemma A.3 of \cite{BuLiv}), we get
$$
|\xi - \alpha| \ll_n {|P(\xi)| \over H(P)} \ll_n H(\alpha)^{-1} \,
(H^{w-n+1})^{-w/(w-n+1)}
$$
and
$$
H(\alpha) \ll_n H(P) \ll_n H^{w-n+1}.
$$
Since all of this is true for every sufficiently large $H$, we get 
$\hw_n^*(\xi) \ge w/(w-n+1)$.  
Selecting now $\eps$ arbitrarily close to $0$, we obtain (2.3).  

In order to establish (2.4), we may assume that
$\hw_n (\xi)$ is finite and set $w = {\hat w}_n(\xi) (1 + \eps)^2$. 
We follow the same argument as in the proof of (2.3).  
The definition of  ${\hat w}_n(\xi)$ and the first inequality of (2.5) then show 
that there exist arbitrarily large values of $H$
for which the polynomial $P(X)$ satisfies $H(P) > H^{1 + \eps}$. We conclude
that there exist algebraic numbers $\alpha$ of arbitrarily large height with
$$
|\xi - \alpha| \ll_n H(\alpha)^{-1 -w/(w-n+1)}.
$$
Thus, we get $w_n^*(\xi) \ge w/(w-n+1)$ and,
selecting $\eps$ arbitrarily close to $0$, we obtain (2.4).  
\cqfd

The next result shows that if $w_n (\xi) = n$ holds, then the values of the five
other exponents at the point $\xi$ are known.

\proclaim Corollary 2.7. 
Let $n$ be a positive integer and $\xi$ a real number such that
$w_n (\xi) = n$. Then we have
$$
w_n(\xi) = w_n^*(\xi) =   \hat{w}_n (\xi) =  
\hat{w}_n^* (\xi) = n    \eqno (2.6)
$$
and
$$
\lambda_n (\xi) = \hat{\lambda}_n (\xi) = {1 \over  n}.    \eqno (2.7)
$$

\pro 
Equalities (2.6) follow from Theorems 2.5 and 2.6. 
Khintchine's transference theorem (see Theorem 3.2 below) shows that
$w_n (\xi) = n$ is equivalent to $\lambda_n (\xi) = 1/n$. Combined
with Theorem 2.5, this gives (2.7). \cqfd

Since $w_n (\xi)$ is equal to $n$ for almost all real numbers $\xi$,
with respect to the Lebesgue measure
(this was proved by Sprind\v zuk \cite{SprLiv}), the next result is
an immediate consequence of Corollary 2.7.

\proclaim Theorem 2.8. 
For almost all (with respect to Lebesgue measure) real numbers
$\xi$ and every positive integer $n$, we have
$$
w_n(\xi) = w_n^*(\xi) =   \hat{w}_n (\xi) =  \hat{w}_n^* (\xi) = n.
$$
and
$$
\lambda_n (\xi) = \hat{\lambda}_n (\xi) = {1 \over  n}.
$$

Before discussing the values taken by our exponents at algebraic points, we recall
Liouville's inequality (see Theorem A.1 and Corollary A.2 in \cite{BuLiv}).

\proclaim Theorem 2.9 (Liouville's inequality). 
Let $\alpha$ and $\beta$ be distinct 
real algebraic numbers of degree $n$ and $m$, respectively. Then,
$$
|\alpha - \beta| \ge (n+1)^{-m} (m+1)^{-n} H(\alpha)^{-m} H(\beta)^{-n}.
$$
Furthermore, if $P(X)$ is an integer polynomial of degree $n$ which does not
vanish at $\beta$, then
$$
|P( \beta)| \ge (n+1)^{-m} (m+1)^{-n}  H(P)^{-m+1} H(\beta)^{-n}.
$$

Let $\xi$ be a real algebraic number of degree $d \ge 1$. 
It follows from Theorem 2.9 that $w_n^* (\xi) \le w_n(\xi) \le d-1$ holds for $n \ge 1$. 
Roth's theorem, which asserts that $w_1 (\xi) = 1$ if $\xi$
is irrational, has been considerably 
extended by Schmidt \cite{SchmLN}, who showed that $w_n (\xi) = w_n^* (\xi) = n$
holds for $n \le d-1$. This enables us to get the values of our six exponents at real algebraic 
numbers.

\proclaim Theorem 2.10. 
Let $\xi$ be a real algebraic number
 of degree $d \ge 1$ and let $n$ be a positive integer. We have
$$
w_n(\xi) = w_n^*(\xi) = \hat{w}_n (\xi) =  \hat{w}_n^* (\xi) = \min\{n, d-1\}
$$
and
$$
\lambda_n (\xi) = \hat{\lambda}_n (\xi) = {1 \over  \min\{n, d-1\}}.
$$

Theorem 2.10 shows that 
real algebraic numbers of degree greater than $n$ do behave like almost all real
numbers, as far as approximation by algebraic numbers of degree
less than $n$ is concerned. 
We may as well consider approximation to complex (non-real) numbers.
Quite surprisingly, complex non-real numbers of degree greater than $n$ do 
not always behave like almost all complex numbers; see \cite{BuEv09}.

Theorem 2.10 shows that we can focus on the values taken by
our exponents at transcendental, real numbers. 
This motivates the following definition.

\proclaim Definition 2.11. 
The spectrum of an exponent of approximation
is the set of values taken by this exponent at transcendental
real numbers.

We point out an important problem, which will be
discussed in the next sections. 

\proclaim Problem 2.12.
To determine the spectra of the exponents 
$w_n, w_n^*, \lambda_n, \hw_n, \hw_n^*, \hla_n$.

Results towards Problem 2.12 are summarized at the end of Section 5.

\vskip 5mm

\centerline{\bf 3. Relations between exponents}

\vskip 5mm

We begin this section with an easy result on the difference between the exponents 
$w_n$ and $w^*_n$. 

\proclaim Theorem 3.1.
For any positive integer $n$ and any transcendental real number $\xi$, we have
$$
w_n(\xi) - n +1 \le w_n^* (\xi) \le w_n (\xi)   \eqno (3.1)
$$
and
$$
\hw_n(\xi) - n +1 \le \hw_n^* (\xi) \le \hw_n (\xi).    \eqno (3.2)
$$

\pro
The right hand side inequalities of (3.1) and (3.2) 
have been already stated in Theorem 2.5. 
The left hand side of inequality (3.1) is inequality (3.11) in \cite{BuLiv}, whose 
proof also gives the left hand side of inequality (3.2). 
\cqfd

It is interesting to note that the left hand side inequality of (3.2) 
is sharp since there exist real numbers $\xi$ with
$\hw_n(\xi) = n$ and $\hw_n^* (\xi) = 1$; see Corollary 5.4 below.  
We do not know if the left hand side inequality of (3.1) is sharp for $n \ge 4$; see
Theorem 5.7 below.

We indicate now some transference results linking together the
  rational simultaneous approximation to $\xi, \ldots, \xi^n$
and the smallness  of the linear form  with coefficients $\xi,
\ldots, \xi^n$.

\proclaim Theorem 3.2.
For every integer $n \ge 2$ and every real number $\xi$ which is not
algebraic of  degree  $\le n$, we have
$$
{1 \over n} \le {w_n (\xi) \over (n-1) w_n (\xi)  + n} \le \lambda_n (\xi) \le
{w_n (\xi) - n  +1 \over n}  \eqno (3.3)
$$
and
%$$
%{1 \over n} \le {\hw_n (\xi) \over (n-1) \hw_n (\xi)  + n} \le \hlambda_n (\xi) \le
%{\hw_n (\xi) - n  +1 \over n}, \eqno (3.5)
%$$
$$
{1 \over n} \le {\hw_n (\xi) - 1 \over (n-1) \hw_n (\xi)} \le \hlambda_n (\xi) \le
{\hw_n (\xi) - n  +1 \over   \hw_n (\xi) }.  \eqno (3.4)
$$

\proof 
The inequalities (3.3) follow direcly from Khintchine's 
transference principle (cf. Theorem B.5 from \cite{BuLiv}),
whose proof shows that the same inequalities hold for the uniform exponents; see \cite{Jar38}. 
The latter result is weaker than inequalities (3.4), which have been
recently established by German \cite{Ger12a}. 
\cqfd

Observe that (3.4) with $n=2$ reduces to the following 
result established by Jarn\'\i k \cite{Jar38} in 1938; see also \cite{Kh48b} for an
alternative proof.

\proclaim Theorem 3.3.
For every transcendental real number $\xi$ we have
$$
\hla_2 (\xi) = 1 - {1 \over \hw_2 (\xi)}.
$$

Inequalities (3.3) have been recently refined in \cite{Lau09b,BuLa10} by means of the
introduction of uniform exponents.

\proclaim  Theorem 3.4. 
Suppose $n\ge 2$. The inequalities 
$$
 \lambda_n(\xi) \ge 
{ (\hw_n(\xi) -1) w_n(\xi)  \over 
((n-2)\hw_n(\xi) +1)w_n (\xi) +(n-1)\hw_n(\xi)}
$$
and
$$
\lambda_n(\xi) \le {(1 - \hlambda_n(\xi))w_n (\xi) -n +2 -\hlambda_n(\xi) \over n-1}
$$
hold for any transcendental real number $\xi$. 

Since $\hla_n(\xi) \ge 1/n$ and $\hw_n (\xi) \ge n$, one easily checks that
Theorem 3.4 implies (3.3).

The first inequality in the next theorem was established by 
Davenport and Schmidt \cite{DaSc69}, while investigating 
the approximation to a real number by algebraic integers.
The second one is a recent result of Schleischitz \cite{Schl14}.

\proclaim Theorem 3.5.
For any positive integer $n$ and any transcendental real number $\xi$,
we have
$$
w_n^* (\xi) \ge {1 \over \hla_n (\xi)}
\quad
\hbox{and}
\quad
\hw_n^* (\xi) \ge {1 \over \lambda_n (\xi)}.
$$

It is interesting to note that the combination of (3.4) with the first inequality of Theorem 3.5 
gives (2.4). 

Theorem 3.5  relates
approximation to $\xi$ by algebraic numbers of degree at most $n$ with
uniform simultaneous rational approximation to $\xi, \ldots, \xi^n$.
Further explanations can be found in \cite{BuTeu} and in Section 3.6 of \cite{BuLiv}.

Note that Theorem 3.5 can be compared with the next result, extracted from 
\cite{BuLa05b,Cas}, which is a particular case of Theorem 1.3. 
If $A$ is the matrix $(\xi, \xi^2, \ldots , \xi^n)$ and ${{\underline \theta}}$ is
a real $n$-tuple, we simply write $w_n  (\xi, {{\underline \theta}})$
for $w  (A, {{\underline \theta}})$ and $\hat{w}_n  (\xi, {{\underline \theta}})$
for $\hat{w}  (A, {{\underline \theta}})$. 

\proclaim Theorem 3.6.
For any integer $n \ge 1$, any transcendental real number $\xi$ and any 
real $n$-tuple ${{\underline \theta}}$, we have
$$
w_n  (\xi, {{\underline \theta}}) \ge {1 \over \hla_n (\xi)}
\quad
\hbox{and}
\quad
\hw_n  (\xi, {{\underline \theta}}) \ge {1 \over \lambda_n (\xi)}.
$$

We continue with a recent result of Marnat \cite{Mar15}
obtained by using the deep result of Roy \cite{Roy15a} on the 
Schmidt--Summerer parametric geometry of numbers; see Section 7.

\proclaim Theorem 3.7.
For $n \ge 3$, there does not exist any algebraic relation between the exponents
$\hw_n$ and $\hla_n$.

The next results were proved by Schmidt and Summerer \cite{ScSu13a,ScSu13b}; see
also Moshchevitin \cite{Mosh10,Mosh12}. 

\proclaim Theorem 3.8. 
For any transcendental real number $\xi$ we have
$$
w_n (\xi) \ge \hw_n(\xi) \, {(n-1) (\hw_n(\xi) - 1) \over 1 + (n-2) \hw_n(\xi)}, 
$$
for $n \ge 2$, and 
$$
w_3 (\xi) \ge \hw_3(\xi) \, {\sqrt{4 \hw_3(\xi) - 3} - 1  \over 2}.
$$

The case $n=2$ of Theorem 3.8 was proved by Jarn\'\i k \cite{Jar54}. 

\proclaim Theorem 3.9.
For any transcendental real number $\xi$ we have
$$
\lambda_n (\xi) \ge \hlambda_n (\xi) \, {\hla_n(\xi)  + n - 2 \over (n-1) (1 - \hla_n(\xi))}, 
$$
for $n \ge 2$, and 
$$
\lambda_3 (\xi) \ge  \hlambda_3 (\xi) \, {\hla_3(\xi)  + \sqrt{\hla_3(\xi) (4 - 3 \hla_3(\xi))} 
\over 2 (1 - \hla_3(\xi))}. 
$$

Actually, Theorems 3.2 to 3.4 and 3.6 to 3.9 
are valid for general tuples $\uxi$, whose coordinates
are, together with $1$, linearly independent over $\Z$ (and not only
for tuples of the shape $(\xi, \xi^2, \ldots , \xi^n)$.)

\vskip 5mm

\centerline{\bf 4. Bounds for the exponents}

\vskip 5mm

Let $n$ be a positive integer. 
As we have seen in Section 2, the Dirichlet {\it Schubfachprinzip} (or, if one prefers,
Minkowski's theorem) readily implies that $w_n(\xi)$ is at least equal to $n$ for
any positive integer $n$ and any real number $\xi$ not algebraic of degree $\le n$.
It is a longstanding problem, which was first formulated by Wirsing \cite{Wir}, 
to decide whether the same
result remains true for the quantity $w_n^*(\xi)$.

\proclaim Conjecture 4.1 (Wirsing). 
For any positive integer $n$ 
and any transcendental real number $\xi$ we have $w_n^*(\xi) \ge n$.

The seminal paper of Wirsing \cite{Wir} and the study of his
conjecture, which has been up to now confirmed only for $n=1$ 
(this follows from the theory of continued fractions) 
and $n=2$ (by Davenport and Schmidt \cite{DaSc67}),
have motivated many works.

\proclaim Theorem 4.2. 
For any real number $\xi$ which is neither rational,
nor quadratic, and for any real number $c$ greater
than $160/9$, there exist infinitely many 
rational or quadratic real numbers $\al$ satisfying
$$
|\xi - \alpha| \le c \, \max\{1, |\xi|^2\} \, H(\alpha)^{-3}.
$$

Theorem 4.2 was proved in \cite{DaSc67} (see also \cite{SchmLN}). 
It has been extended by Davenport and Schmidt \cite{DaSc68}  
(up to the value of the numerical constant) as follows.

\proclaim Theorem 4.3. 
Let $n \ge 2$ be an integer 
and let $\xi$ be a real number which is not algebraic of degree at most $n$. 
Then there exist an effectively computable
constant $c$, depending only on $\xi$ and on $n$, an integer $d$
with $1 \le d \le n-1$, and infinitely many integer polynomials
$P(X)$ of degree $n$ whose roots $\alpha_1, \ldots,
\alpha_n$ can be numbered in such a way that
$$
|(\xi- \alpha_1) \ldots (\xi - \alpha_d)| \le c \, H(P)^{-n-1}.
$$

Theorem 4.2 has recently been improved by Moshchevitin \cite{Mosh14} as follows.

\proclaim Theorem 4.4.
For any real number $\xi$ which is neither rational, nor a quadratic irrationality,
we have
$$
w_2^* (\xi) \ge \hw_2 (\xi) (\hw_2 (\xi) - 1) \ge 2. 
$$

The proof of Theorem 4.4 combines ideas from \cite{DaSc67}
with an argument used by Jarn\'\i k \cite{Jar50,Jar54}
in his proof of the case $n=2$ of Theorem 3.8.

%an argument of Jarn\'\i k \cite{Jar50,Jar54}, who established that,
%if $\xi_1, \xi_2$ are linearly independent real numbers, then we have
%$$
%w_2 (\xi_1, \xi_2) \ge \hw_2 (\xi_1, \xi_2) (\hw_2 (\xi_1, \xi_2) - 1).    \eqno (4.4)
%$$

The first statement of Theorem 4.5 was proved by Wirsing \cite{Wir} and the
second one by Bernik and Tishchenko \cite{BeTi}. 

\proclaim Theorem 4.5.
Let $n$ be a positive integer and $\xi$
be a transcendental real number. Then we have
$$
w_n^*(\xi) \ge {w_n(\xi) + 1 \over 2},  
$$
and
$$
w_n^*(\xi) \ge {n \over 4} + {\sqrt{n^2 + 16 n - 8} \over 4}. \eqno (4.1)
$$

For proofs of Theorem 4.5 and related results, 
the reader may consult Chapter 3 of \cite{BuLiv}.
Slight improvements on (4.1) have been subsequently obtained by 
Tishchenko \cite{Tish00,Tish07}, with very technical proofs.

Davenport and Schmidt \cite{DaSc69} have given uniform upper bounds 
for the exponents $\hla_n$ and $\hw_n$. Some of their results have 
been subsequently  improved by Laurent \cite{Lau03}. 
For a positive real number $x$, we denote by 
$\lceil x \rceil $ the smallest integer greater than or equal to $x$. 
The next theorem expresses results of \cite{DaSc69,Lau03} in terms of exponents 
of approximation. 

\proclaim Theorem 4.6. 
For any positive integer $n$ and any transcendental
real number $\xi$, we have
$$
\hla_n (\xi) \le {1 \over  \lceil n/2  \rceil} 
\quad \hbox{and} \quad \hw_n^* (\xi) \le \hat{w}_n (\xi) \le 2n - 1.
$$

\pro
We only explain how to get the easy estimate $\hw_n^* (\xi) \le 2n - 1$. 
Let $w< \hw_n^*(\xi)$ be a real number. By definition of $\hw_n^*$, there exist 
arbitrarily large integers $H$ such that there are distinct algebraic numbers
$\alpha_1, \alpha_2$ of degree at most $n$ and height at most $H$
such that
$$
|\xi - \alpha_1| < H(\alpha_1)^{-1} \, H^{-w}, 
\quad |\xi - \alpha_2| < H(\alpha_2)^{-1} \, H^{-w}.
$$
This implies that $|\alpha_1 - \alpha_2| < 2 H^{-w-1}$, while 
Theorem 2.9 (Liouville's inequality) gives
that $|\alpha_1 - \alpha_2| \gg_n H^{-2n}$. By combining these two inequalities,
we get that $w \le 2n - 1$. This proves the upper bound $\hw^*_n (\xi) \le 2n - 1$. 
Arguing now with polynomials instead of algebraic numbers, this can be strengthened 
to $\hw_n (\xi) \le 2n - 1$, as was shown in \cite{DaSc69}. 
\cqfd

\proclaim Theorem 4.7.
For any transcendental real number $\xi$, we have
$$
\hla_2 (\xi) \le {\sqrt{5} - 1 \over 2} 
\quad \hbox{and} \quad  \hw_2 (\xi) \le {3 + \sqrt{5}  \over 2},   \eqno (4.2)
$$
and both inequalities are best possible.  
For any transcendental real number $\xi$, we have
$$
\hla_3(\xi) \le \Bigl(2 + \sqrt{5} - \sqrt{7 + 2 \sqrt{5}} \, \Bigr)/2 = 0.4245 \ldots 
$$

The bound for $\hla_2 (\xi)$ given in (4.2) was proved in \cite{DaSc69} and the one
for $\hw_2$ in \cite{ArRo}, whose
authors were at that time not aware of Jarn\'\i k's result (Theorem 3.3).

The last assertion of Theorem 4.7 was proved by Roy \cite{Roy08b}, who improved the
estimate $\hla_3(\xi) \le 1/2$ given in Theorem 4.6. 
He further indicated that his upper bound for the exponent is not best possible.

Roy \cite{Roy03a,Roy04a} showed that the inequalities (4.2) are sharp. 
The set of values taken by the exponents $\hw_2$ and $\hw_2^*$ 
has been studied in \cite{BuLa05a,Fisc04,Fisc06,Fisc07,Roy04b,Roy07,Roy08a}. 
Among other results,
we know that the spectrum of $\hw_2$ is dense in $[2, (3 + \sqrt{5})/2]$ and that
there exists a real number $c < (3 + \sqrt{5})/2$ such that the intersection 
of this spectrum with $[c, (3 + \sqrt{5})/2]$ is countable.

The next result provides, under a suitable assumption, 
an upper bound for $\hw_n^*$ in terms of $w_m$, when $m$ is less than $n$. 
It has not been noticed previously.

\proclaim Theorem 4.8. 
Let $m, n$ be positive integers with
$1 \le m \le n-1$ and $\xi$ be a real number.
If $w_m (\xi) \ge m + n - 1$, then we have
$$
\hw_n^* (\xi) \le {m w_m (\xi) \over w_m (\xi) - n +1}.
$$

%The proof of Theorem 4.8 is inspired by that of Proposition 2.1 of \cite{BuLa05a}. 
%\noindent {\it Proof of Theorem 4.8.} 

\pro 
It is inspired by that of Proposition 2.1 of \cite{BuLa05a}. 
Let $m, n$ and $\xi$ be as in the statement of the theorem. 
Assume for convenience that $|\xi| \le 1$. 
Let $\eps$ be a real number with $0 < \eps \le 1/2$. 
Set $w := w_m (\xi)$.
Let $P(X)$ be an integer polynomial of degree at most $m$
and height $H := H(P)$ large such that  
$$
H^{-w-\eps} < |P(\xi)| < H^{-w + \eps}.    \eqno (4.3)
$$
By using Lemma A.3 from \cite{BuLiv}, we may assume without any loss of
generality that $P(X)$ is irreducible and primitive. 
Let $v$ be a positive real number and set $X = H^{w / (v(1 + \eps))}$. 
%to be solved in algebraic numbers $\alpha$ with $H(\alpha) \le X$ and degree
%at most $n$. 
Let $\alpha$ be the root of $P(X)$ which is the closest to $\xi$. 
If $|\xi - \alpha| \le H(\alpha)^{-1}  H^{- w(1 + 2 \eps)/(1 + \eps)}$, then,
by Corollary A.1 of \cite{BuLiv}, we have
$$
|P(\xi)| \ll_n |\xi - \alpha| \cdot H(P) \ll_n H^{-w (1 + 2 \eps/3)} \ll_n H^{-w - 4\eps/3},
$$
using that $\eps \le 1/2$ and $w \ge 2$. This gives a 
contradiction to (4.3) if $H$ is large enough. Consequently, we have
$$
|\xi - \alpha| > H(\alpha)^{-1}  X^{- v(1 + 2 \eps)},     \eqno (4.4)
$$ 
provided that $H$ is large enough. 

Let $\beta$ be an algebraic number of height at most $X$ 
and degree at most $n$, which satisfies 
$$
|\xi - \beta| \le H(\beta)^{-1} X^{-v(1+ 2 \eps)}.   \eqno (4.5) 
$$
It follows from (4.4) that $\beta \not= \alpha$.
Liouville's inequality  (Theorem 2.9) then gives
$$
|P(\beta)| \gg_n  H^{-n+1} \, H(\beta)^{-m}.   \eqno (4.6)
$$ 
By Rolle's theorem and the fact that $|\xi| \le 1$, we have
$$
|P(\beta)| \le |\xi - \beta| \cdot \max_{t : |t| \le 2} |P'(t)| + |P(\xi)| \le 
m 2^m |\xi - \beta|  H + H^{-w+\eps}.     \eqno (4.7) 
$$
If $H^{-w+\eps} \ge |\xi - \beta| \cdot H$, then (4.6) and (4.7) imply  
$$
H(\beta)^{-m} \ll_n H^{n-1-w+\eps}
$$
and, since $H(\beta) \le H^{w/(v(1+\eps))}$, we get
$$
v \le {mw \over (w+1-n - 2\eps) (1 + \eps)},  \eqno (4.8)
$$
provided that $H$ is large enough. 

If $H^{-w+\eps} \le |\xi - \beta| \cdot H$, then, by (4.5), (4.6) and (4.7), we get
$$
H(\beta)^{-m+1} H^{-n} \ll_n H^{-w(1 + 2 \eps)/(1 + \eps)},
$$
hence,
$$
H^{(m-1) w / (v(1+\eps))} \gg_n H^{-n + w(1 + 2 \eps)/(1 + \eps)}.
$$
This implies
$$
v \le {(m-1)w \over w  (1+ 2 \eps) - n   (1 + \eps) - \eps}, \eqno (4.9)
$$
provided that $H$ is large enough. 
Consequently, we have $\hw_n^* (\xi) < v (1 + 2 \eps)$ as soon as 
$v$ exceeds the left-hand sides of (4.8) and (4.9).  
Since $\eps$ can be taken arbitrarily close to $0$, we deduce that
$$
\hw_n^* (\xi) \le \max \Bigl\{ {mw \over w + 1 - n}, {(m-1)w \over w-n} \Bigr\}.
$$
Since $w \ge m + n - 1$, this maximum is equal to $mw / (w+1-n)$.
This ends the proof of the theorem. 
\cqfd

We end this section with a recent result of Schleischitz \cite{Schl15a}.

\proclaim Theorem 4.9.
For any positive integer $n$ and any transcendental real number $\xi$, we have
$$
\hlambda_n (\xi) \le \max\Bigl\{ {1 \over n}, {1 \over \lambda_1 (\xi)} \Bigr\}.
$$

\pro
Assume that $n \ge 2$ and that $\xi$ is in $(0, 1)$ with $\lambda_1 (\xi) > 1$.  
Let $q$ be a large positive integer and $v$ be a real number greater than $1$
such that $q^{-v} < 1/(4 n q)$ and
$$
\| q \xi \|  \le q^{-v}.        \eqno (4.10)
$$
Let $p$ be the integer such that $|q \xi - p| = \| q \xi \|$.
Without any loss of generality, we may assume that $p$ and $q$ are coprime. 
Observe that, for $j = 1, \ldots , n$, we have
$$
\Bigl| \xi^j - {p^j \over q^j} \Bigr| 
=  \Bigl| \xi  - {p  \over q} \Bigr| \cdot
\Bigl| \xi^{j-1} + \ldots + \Bigl(  {p  \over q} \Bigr)^{j-1} \Bigr| 
\le {2n  \over q^{1+v}} \le {1 \over 2 q^2}.  \eqno (4.11)
$$
Let $v'$ be a real number with $1 < v' < \min\{v, n\}$ and set $X = q^{v'}$. 
Let $x$ be a positive integer with $x < X$ and express $x$ in base $q$. 
There exist integers $b_0, b_1, \ldots , b_{n-1}$ in $\{0, 1, \ldots , q - 1\}$
such that
$$
x = b_0 + b_1 q + b_2 q^2 + \ldots + b_{n-1} q^{n-1}.
$$
Let $u$ be in $\{1, 2, \ldots , n\}$ be the smallest index
such that $b_{u-1}$ is non-zero. Then, 
$$
\| x p^u / q^u\| = \|b_{u-1} p^u / q\| \ge 1/q.   \eqno (4.12)
$$
Furthermore, it follows from (4.11) and the fact that $v' < v$ that
$$
\Bigl| x \Bigl(\xi^u -  {p^u  \over q^u} \Bigr) \Bigr| 
\le {2 n q^{v'} \over q^{1 + v}}  < {1 \over 2q},
$$
if $q$ is sufficiently large. Denoting by $y$ the nearest integer to $x \xi^u$, 
the triangle inequality and (4.12) then give 
$$
\eqalign{
\max_{1\le j \le n} \| x \xi^j \| \ge \| x \xi^u \| 
& = \Bigl| \Bigl( {x p^u  \over q^u} - y \Bigr) + x \Bigl(\xi^u -  {p^u  \over q^u} \Bigr)   \Bigr|  \cr
& \ge  {1 \over q} - {1 \over 2 q} 
\ge {1 \over 2 q} =  {1 \over 2 X^{1/v'}}. \cr}
$$
If (4.10) is satisfied for arbitrarily large integers $q$, this shows that
$$
\hlambda_n (\xi) \le {1 \over v'}.
$$
As $v'$ can be chosen arbitrarily close to $\min\{\lambda_1 (\xi), n\}$, we have proved that
$$
\hlambda_n (\xi) \le {1 \over \min\{\lambda_1 (\xi), n\}} = 
\max\Bigl\{ {1 \over n}, {1 \over \lambda_1 (\xi)} \Bigr\}.
$$
This completes the proof of the theorem. \cqfd 

We point out an immediate corollary of Theorem 4.9. 

\proclaim Corollary 4.10.
For any positive integer $n$, any transcendental real number $\xi$
with $w_1 (\xi) \ge n$ satisfies $\hw_k (\xi) = k$ for $k = 1, \ldots , n$.

\pro
This follows from Theorem 4.9 combined with Theorem 3.2. \cqfd

\vskip 5mm

\centerline{\bf 5. Spectra}

\vskip 5mm

This section is mainly devoted to the study of the spectra of the six exponents 
of approximation defined in Section 2. 

We begin with an auxiliary result,
extracted from \cite{Bu10a}, which confirms the existence of real numbers $\xi$
for which $w_1 (\xi) = w_n (\xi)$, for some integer $n \ge 2$. 

\proclaim Theorem 5.1.
Let $n \ge 1$ be an integer.
For any real number $w \ge 2n-1$, there exist uncountably many 
real numbers $\xi$ such that
$$
w_1(\xi) = w_1^*(\xi) = \ldots = w_n (\xi) 
= w_n^*(\xi) = w.
$$
%and
%$$
%\hw_j (\xi) = j, \quad \hbox{for $j=1, \ldots , n$}.
%$$
In particular, the spectra of $w_n$ and $w_n^*$ include the interval $[2n-1, + \infty]$. 

\pro
This is clear for $w=n=1$. 
Let $w>1$ be a real number.
Let $M$ be a large positive integer
and consider the real number
$$
\xi_ w = [0; 2, M \lfloor q_1^{w-1} \rfloor , M \lfloor q_2^{w-1} \rfloor , 
M \lfloor q_3^{w-1} \rfloor , \ldots ],   
$$
where $q_1 = 2$ and $q_j$ is the denominator of the $j$-th
convergent to $\xi_ w$, that is, of
the rational number 
$p_j / q_j = [0; 2, M \lfloor q_1^{w-1} \rfloor , \ldots , 
M \lfloor q_{j-1}^{w-1} \rfloor]$, for $j \ge 2$.
By construction, we have
$$
q_{j+1} \asymp M q_j^w \quad
\hbox{and} \quad
\biggl | \xi_ w - {p_j \over q_j} \biggr | \asymp {1 \over M q_j^{w+1}}, 
$$
for $j \ge 1$. Consequently, we have
$$
w = w_1 (\xi_ w)  \le \ldots \le w_n (\xi_ w). \eqno (5.1)
$$

Using triangle inequalities, it is shown in \cite{Bu10a} that, if $M$
is sufficiently large and $w \ge 2 n - 1$, then
$$
|P (\xi_ w) | \gg_{n, M} H(P)^{-w}
$$
holds for every
polynomial $P(X)$ of degree at most $n$ and sufficiently large height,
hence $w_n (\xi_ w) \le w$ and the inequalities in (5.1) are indeed equalities. 
An additional argument is needed to show that 
$w_1^*(\xi) = \ldots = w_n^*(\xi) = w$; see \cite{Bu10a} for the complete proof. 
\cqfd

It would be desirable to replace the assumption $w \ge 2n-1$
in Theorem 5.1 by a weaker one. Actually, the value $2n-1$ comes from 
Theorem 2.9 (Liouville's inequality),  
which is widely used in the proof of Theorem 5.1.

Theorem 5.1 is a key tool to get results on the spectra of various exponents
of approximation.

The next result, also established in \cite{Bu10a}, gives a relationship between 
the exponents $\lambda_n$ and $\lambda_m$
when $m$ divides $n$.

\proclaim Lemma 5.2.
For any positive integers $k$ and $n$, and any
transcendental real number $\xi$ we have
$$
\lambda_{kn} (\xi) \ge {\lambda_k (\xi) - n + 1 \over n}.
$$

A similar inequality holds between the uniform exponents, but it gives nothing interesting
since $\hla_k (\xi) \le 1$ for every integer $k \ge 1$ and every irrational real
number $\xi$. 

\pro
Let $v$ be a positive real number
and $q$ be a positive integer such that
$$
\max_{1\le j \le k} | q \xi^j - p_j|  \le q^{-v},
$$
for suitable integers $p_1, \ldots , p_k$.
Let $h$ be an integer with $1 \le h \le kn$. 
Write $h = j_1 + \ldots + j_m$ with $m \le n$ and
$1 \le j_1, \ldots , j_m \le k$.
Then,
$$
|q^m \xi^h - p_{j_1} \ldots p_{j_m}| \ll_m q^{m-1} q^{-v}
$$
and
$$
\| q^n \xi^h \| \ll q^{n-m} \| q^m \xi^h \| \ll_m
q^{n-1-v} \ll_m (q^n)^{- (v-n+ 1)/n},
$$
independently of $h$.
This proves the lemma. \cqfd

We display an immediate consequence of Lemma 5.2.

\proclaim Corollary 5.3.
Let $\xi$ be a real irrational number.
Then, $\lambda_n (\xi) = + \infty$ 
holds for every positive integer $n$
if, and only if, $\lambda_1 (\xi) = + \infty$. 

Combined with Theorems 2.5 and 4.8 and Corollary 4.10, Corollary 5.3
allows us to determine 
the values taken at Liouville numbers (recall that a Liouville number is, by definition,
a real number $\xi$ satisfying $w_1 (\xi) = +  \infty$) 
by our six exponents of approximation.

\proclaim Corollary 5.4.
For any positive integer $n$ and any Liouville number $\xi$, we have
$$
w_n(\xi) = {w}_n^* (\xi) = \lambda_n (\xi) = + \infty, \quad \hw_n (\xi) = n, \quad
\hat{w}_n^* (\xi) = 1, \quad
\hbox{and} \quad \hla_n (\xi) = {1 \over n}.
$$

The proof of Theorem 5.1 shows how the theory of continued fractions
allows us to construct explicitly real numbers
$\xi$ having any arbitrarily prescribed value 
for $\lambda_1 (\xi)$.
The same question for an exponent $\lambda_n$ with $n \ge 2$
is not yet solved. At present, the best known result was proved in \cite{Bu10a}
and is reproduced below.

\proclaim Theorem 5.5. 
Let $n \ge 2$ be an integer and $w \ge n$ be a real number.
If the real number $\xi_w$ satisfies
$w_1(\xi_w) = \ldots = w_n(\xi_w) = w$, then
$$
\lambda_n (\xi_w) = {w - n + 1 \over n}, \quad \hw_n^* (\xi_w) = {w \over w - n +1},
$$
and
$$
\hw_j (\xi_w) = j, \quad j= 1, \ldots , n.
$$

\pro
Let $m \ge 2$ be an integer and $\xi$ be a
transcendental real number.
Lemma 5.2 with $k=1$ implies the lower bound
$$
\lambda_m(\xi) \ge { w_1(\xi) -m+1\over m}. 
$$
On the other hand, Theorem 3.2
gives the upper bound 
$$
\lambda_m(\xi) \le { w_m(\xi) -m +1\over m}.  
$$
Let $\xi_w$ be such that
$$
w_1(\xi_w) = \ldots = w_n(\xi_w) = w.  
$$
Then,  the equalities
$$
\lambda_m (\xi_w) = {w - m + 1 \over m},  
\quad m= 1, \ldots , n,
$$
hold; in particular,
$$
\lambda_n (\xi_w) = {w - n + 1 \over n},
$$
and this establishes the first statement of the theorem. 

Combining Theorem 2.6 with the case $m=1$ of Theorem 4.8 gives
$$
{w_n(\xi_w)  \over w_n(\xi_w)  - n +1}  \le \hw_n^* (\xi_w) 
\le  {w_1 (\xi_w)  \over w_1(\xi_w)  - n +1},
$$
thereby proving the second statement of the theorem.

Without any loss of generality, assume that $0 < \xi_w < 1$ and $w > n \ge 2$. 
Let $\eps$ be a real number satisfying $0 < \eps < w-n$.
Let $p/q$ be a rational number
such that $q \ge 1$, $\gcd(p, q) = 1$
and $|\xi_w - p/q| < q^{-1-w + \eps}$. Let $P(X)$ be an integer 
polynomial of degree $j$ at most $n$ and height at most $q-1$. 
Observe that $P(p/q)$ is a non-zero rational number satisfying 
$|P(p/q)| \ge 1/q^j$. 

By Rolle's theorem, there exists $t$ lying between $\xi$ and  $p/q$ such that
$$
P(\xi_w) = P(p/q) + (\xi_w - p/q) P'(t).
$$
Observe that $|P'(t)| \le 2^n n^2 q$ and $|\xi_w - p/q| \cdot |P'(t)| \le 1 / (2 q^n)$,
if $q$ is large enough.
We then deduce that
$|P(\xi_w) | \ge 1 / (2 q^j)$ if $q$ is large enough.
This shows that $\hw_j (\xi_w) = j$, as asserted.
\cqfd

%\proclaim Corollary 5.6.
%Let $n$ be a positive integer and $\lambda \ge 1$ a real number.
%There are uncountably many real numbers $\xi$, which can be
%constructed explicitly, such that $\lambda_n (\xi) = \lambda$.
%In particular, the spectrum
%of $\lambda_n$ includes  the interval $[1, +\infty]$.

%\pro This follows from Theorem 5.5 combined with Theorem 5.1. \cqfd

Theorem 5.5 allows us to get some information on the spectra of the
exponents $\lambda_n$ and $\hw_n^*$.

\proclaim Theorem 5.6.
For any positive integer $n$,  the spectrum of $\lambda_n$ includes  the interval $[1, +\infty]$,
the spectrum of $\hw_n^*$ includes the interval $[1, 2- 1/n]$
and the spectrum of $\hw_n - \hw_n^*$ includes the interval 
$[n-2+1/n, n-1]$.

\pro 
This has been already proved for $n=1$.
If $n \ge 2$, the statement follows from the combination of Theorem 5.1
with Theorem 5.5. 
\cqfd

Recall that, by Theorem 3.1, the spectra of $w_n - w_n^*$  and
of $\hw_n - \hw_n^*$ are included in $[0, n-1]$, for $n \ge 1$. 
The first assertion of the
next result was proved by Bugeaud and Dujella \cite{BuDu11}
by means of an explicit construction of families of polynomials
with close roots.

\proclaim Theorem 5.7.
For any positive integer $n$, the 
spectrum of $w_n - w_n^*$ includes the interval 
$$
\Bigl[0, {n \over 2} + {n-2 \over 4(n-1)} \Bigr).
$$
Moreover, the spectrum of $w_2 - w_2^*$ is equal to $[0, 1]$
and that of $w_3 - w_3^*$ is equal to $[0, 2]$.

Explicit examples of real numbers $\xi$ for which
$w_2 (\xi)$ exceeds $w_2^* (\xi)$ can be found in \cite{Bu12b}.

Very recently, Schleischitz \cite{Schl15a} established that, under some
extra assumption, the inequality proved in Lemma 5.2 is indeed an equality. 

\proclaim Theorem 5.8. 
Let $n$ be a positive integer and $\xi$ be a real number. If $\lambda_n (\xi) > 1$,
then we have
$$
\lambda_1 (\xi) = n \lambda_n (\xi) + n - 1    \eqno (5.2)
$$
and
$$
\lambda_j (\xi) = {n \lambda_n (\xi) - j + n  \over j} , \quad
\hlambda_j (\xi) = {1  \over j} , \quad (j = 1, \ldots , n).
$$
Conversely, if $\lambda_1 (\xi) > 2n - 1$, then we have
$$
\lambda_n (\xi) = {\lambda_1 (\xi) - n + 1 \over n}    \eqno (5.3)
$$
and
$$
\hlambda_j (\xi) = {1 \over j} , \quad (j = 1, \ldots , n).    \eqno (5.4)
$$

\pro
Assume that $\xi$ is in $(0, 1)$ and satisfies $\lambda_n (\xi) > 1$. 
Let $q$ be a large positive integer and $v$ be a real number greater than $1$
such that $q^{-v} < 1/(4nq)$ and
$$
\max_{1\le j \le n} \| q \xi^j \|  \le q^{-v}.  \eqno (5.5)
$$
Let $p$ be the integer such that $|q \xi - p| = \|q \xi \|$. Then, $p$ and $q$ may not be 
coprime, but $p/q$ is a convergent to $\xi$. Let $d$ be the greatest prime factor of
$p$ and $q$ and set $p_0 = p/d$ and $q_0 = q/d$. 
Observe that, for $j = 1, \ldots , n$, we have
$$
\Bigl| \xi^j - {p^j \over q^j} \Bigr| 
=  \Bigl| \xi  - {p  \over q} \Bigr| \cdot
\Bigl| \xi^{j-1} + \ldots + \Bigl(  {p  \over q} \Bigr)^{j-1} \Bigr| 
\le {2 n \over q^{1+v}} \le {1 \over 2 q^2}.  \eqno (5.6)
$$
Assume that $q < q_0^n$ and express $q$ in base $q_0$. Recalling that $q_0$
divides $q$, there exist integers $b_1, \ldots , b_{n-1}$ in $\{0, 1, \ldots , q_0 - 1\}$
such that
$$
q = b_1 q_0 + b_2 q_0^2 + \ldots + b_{n-1} q_0^{n-1}.
$$
Let $u$ be in $\{2, \ldots , n\}$ be the smallest index
such that $b_{u-1}$ is non-zero. Then, 
$$
\| q p_0^u /q_0^u\| = \|b_{u-1} p_0^u / q_0\| \ge 1/q_0.   \eqno (5.7)
$$
Furthermore, it follows from (5.6) that
$$
\Bigl| q \Bigl(\xi^u -  {p^u  \over q^u} \Bigr) \Bigr| 
= \Bigl| q \Bigl(\xi^u -  {p_0^u  \over q_0^u} \Bigr) \Bigr| \le
{1 \over 2 q}.    \eqno (5.8) 
$$
Let $y$ be the integer such that $|q \xi^u - y| = \| q \xi^u\|$ and observe that
$$
\| q \xi^u\| \ge \Bigl| y - q {p_0^u \over q_0^u} \Bigr| 
-  \Bigl|  q \Bigl( \xi^u -  {p_0^u \over q_0^u} \Bigr) \Bigr|,
$$
using the triangle inequality. Combined with (5.7) and (5.8), this gives 
$$
\max_{1\le j \le n} \| q \xi^j \| \ge \| q \xi^u \| \ge  {1 \over q_0} - {1 \over 2 q} 
\ge {1 \over 2 q},
$$
a contradiction to (5.5).  

Consequently, $b_1 = \ldots = b_{n-1} = 0$
and we have established that $q \ge q_0^n$ (actually, our proof
shows that $q$ must be an integer multiple of $q_0^n$). 
In particular, we have
$$
d \ge q_0^{n-1}.   \eqno (5.9) 
$$
Since
$$
\| q_0 \xi \| = |q_0 \xi - p_0| = d^{-1} \| q \xi \|,
$$
it follows from (5.5) and (5.9) that
$$
\|q_0 \xi\|  \le q_0^{-n+1} q_0^{-n v} = q_0^{-n v - n + 1}.
$$
Since $v$ can be taken arbitrarily close to $\lambda_n (\xi)$, we deduce that
$$
\lambda_1 (\xi) \ge n \lambda_n (\xi) + n - 1.
$$
Combined with Lemma 5.2, this proves the first statement of the theorem. 

In particular, we get $\lambda_1 (\xi) > 2n-1$, and it follows from
Theorem 4.9 that $\hla_j (\xi) = 1/j$ for $j = 1, \ldots , n$. 

Let $j$ be an integer with $2 \le j \le n-1$. Since $\lambda_n (\xi) > 1$, we have 
$\lambda_j (\xi) > 1$ and $\lambda_1 (\xi) = j \lambda_j (\xi) + j - 1$. Combined with (5.2),
this gives $j \lambda_j (\xi) = n \lambda_n (\xi) - j + n$, as claimed.

If the real number $\xi$ satisfies $\lambda_1 (\xi) > 2n - 1$, then we get by Lemma 5.2 that
$\lambda_n (\xi) > 1$ and (5.3) and (5.4) follow from the first assertions of the theorem.
\cqfd

The condition $\lambda_n (\xi) > 1$ in the statement of Theorem 5.8 cannot be removed
in view of Theorem 4.3 of \cite{Bu10a}, which confirms 
the existence of uncountably many real numbers
$\xi$ satisfying $\lambda_n (\xi) = 1$ for every $n \ge 1$. 
Furthermore, Theorem 4.4 of \cite{Bu10a} asserts that, for an arbitrary real number 
$\lambda$ in $[1, 3]$, there exist uncountably many  real numbers $\xi$
satisfying $\lambda_1 (\xi) = \lambda$ and $\lambda_2 (\xi) = 1$.

We display an immediate consequence of Theorem 5.8.

\proclaim Corollary 5.9.
Let $n$ be a positive integer and $\xi$ be a transcendental real number.
Then, $\lambda_n (\xi) > 1$ holds if and only if $\lambda_1 (\xi) > 2n - 1$ holds.

The restriction $w \ge 2n-1$ in the statement of Theorem 5.1 prevents us to get the
whole spectra of the exponents $w_n$ and $\lambda_n$ by the method described above.
Actually, we need the help of metric number theory to determine the whole spectra
of the exponents $w_n$ and $\lambda_2$. 

\proclaim Theorem 5.10.
For any positive integer $n$, the 
spectrum of $w_n$ is equal to the whole interval $[n, + \infty]$
and the spectrum of $w_n^*$ includes the whole interval $[n, + \infty]$.

The first statement of Theorem 5.10 was proved by Bernik \cite{Ber83}
and the second one is a result of Baker and Schmidt \cite{BaSc70}.

We display an immediate consequence of results
by Beresnevich, Dickinson, Vaughan
and Velani \cite{BeDiVe07,VaVe06}. 

\proclaim Theorem 5.11. 
The spectrum of $\lambda_2$ is equal to $[1/2, + \infty]$.

More is known than the mere statement of Theorems 5.10 and 5.11.
Indeed, for an integer $n \ge 1$ and a real number $w \ge n$, 
the Hausdorff dimension of the set of real numbers $\xi$ for which 
$w_n(\xi) = w$ (resp., $w_n^* (\xi) = w$) is equal to $(n+1)/(w+1)$.
Furthermore, the Hausdorff dimension of the set of real numbers $\xi$ for which 
$\lambda_2 (\xi) = \lambda$ is equal to $1 / (1 + \lambda)$ if $\lambda \ge 1$ and to
$(2 - \lambda) / (1 + \lambda)$ if $1/2 \le \lambda \le 1$.

The spectra of the exponents $\hw_n$ and $\hla_n$ remain very mysterious
for $n \ge 3$, since we can not even exclude that they are, respectively, reduced
to $\{n\}$ and $\{1/n\}$ (recall that, by spectrum, we mean the set 
of values taken at transcendental points). 

For $n = 2$, the situation is slightly better. By Jarn\'\i k's Theorem 3.3, the value of $\hw_2$
determines that of $\hla_2$, thus it is sufficient to determine the range of
$\hw_2$; see also below Theorem 4.7.  
%Partial results can be found in \cite{BuLa05a}

As for the exponent $\hw_n^*$, it is likely that its spectrum includes the interval
$[1, n]$, but this is not yet proved. 

\vskip 5mm

\centerline{\bf 6. Intermediate exponents}

\vskip 5mm

Let $n \ge 2$ be an integer and $\utheta$ be a point in ${\bf R}^n$.
In \cite{Lau09b}, Laurent introduced new exponents
$\om_{n,d} (\utheta)$ (simply denoted by $\om_d (\utheta)$ in \cite{Lau09b},
since $n$ is fixed throughout that paper) 
measuring the sharpness of the
approximation to $\utheta$ by linear rational varieties
of dimension $d$. He split the Khintchine transference principle
into $n-1$ intermediate estimates 
which connect the exponents $\om_{n,d} (\utheta)$ 
for $d=0, 1, \ldots , n-1$ (see also \cite{BuLa10}).
Actually, Schmidt \cite{Schm67} was the first to investigate
the properties of these exponents $\om_{n,d}$, but he did not
introduce them explicitly.
We briefly recall their definition and we consider new
exponents $w_{n, d}$ defined over ${\bf R}$ by restricting
$\omega_{n,d}$ to the Veronese curve $(x, x^2, \ldots , x^n)$.
It is convenient to view ${\bf R}^n$ as a subset of $\bP^n({\bf R})$
via the usual embedding $(x_1,\dots ,x_n)\mapsto (1,x_1,\dots ,x_n)$. 
We shall  identify $\utheta  =(\xi_1, \dots , \xi_n)$
with its image in $\bP^n({\bf R})$.
Denote by d the projective distance 
on $\bP^n({\bf R})$ and, for any real linear subvariety $L$ 
of $\bP^n({\bf R})$, set
$$
{\rm d} ( \utheta , L) = \min_{ P\in L} {\rm d} (\utheta , P)
$$
the minimal distance between $\utheta$ and the real  points $P$ of $L$. 
When $L$ is rational over ${\bf Q}$, we indicate  
by $H(L)$ its height, that is the Weil height 
of any system of Pl\" ucker coordinates 
of $L$. 
We refer to \cite{Lau09b,BuLa10}
for precise definitions of the projective distance, heights, etc.

\proclaim Definition 6.1.  
Let $n \ge 2$ and $d$ be integers with $0\le d \le n-1$. 
Let $\utheta$ be in ${\bf R}^n$.
We denote by  $\omega_{n,d}(\utheta)$   the supremum  
of the real numbers $\om$ for which there exist 
infinitely many rational linear subvarieties $L \subset \bP^n({\bf R})$
of dimension $d$ such that 
$$
H(L) {\rm d}(\utheta , L) \le H(L)^{-\omega}.
$$
We denote by  $\homega_{n,d}(\utheta)$   the supremum  
of the real numbers $\homega$ such that, for every sufficiently large value of $H$, there exists 
a rational linear subvariety $L \subset \bP^n({\bf R})$
of dimension $d$ with 
$$
H(L) {\rm d}(\utheta , L) \le H^{-\homega}.
$$
If there exists $\xi$ such that $\utheta = (\xi, \xi^2, \ldots , \xi^n)$,
then we set $w_{n, d} (\xi) = \omega_{n,d}(\utheta)$ and 
$\hw_{n, d} (\xi) = \homega_{n,d}(\utheta)$.

We observe that the functions $\lambda_n$ and $w_{n, 0}$ (resp. $w_n$
and $w_{n, n-1}$) coincide. The exponents $\homega_{n,d}$ were 
introduced in \cite{BuLa07,Roy15b}. 

The following transference theorem was proved in \cite{Schm67,Lau09b}. 

\proclaim Theorem 6.2.
Let $n$ be a positive integer.  For any non-zero vector $\uu$ in $\R^{n+1}$,
we have $\omega_{n,0} (\uu)\ge 1/n$ and
$$
{j\omega_{n,j} (\uu) \over \omega_{n,j} (\uu)+j+1}
  \le \omega_{n, j-1}(\uu)
  \le  {(n-j)\omega_{n,j} (\uu)-1 \over n-j+1},
 \quad
j = 1, \ldots ,  n-1,
$$
with the convention that the left-hand side is equal
to $j$ if $\omega_{n,j} (\uu)$ is infinite.

Let the spectrum of the function $\om_{n,d}$ denote the set
of values taken by the exponents $\om_{n,d} (\utheta)$ when $\utheta = (\xi_1, \ldots , \xi_n)$
ranges over ${\bf R}^n$, with $1, \xi_1, \ldots, \xi_n$
linearly independent over the rationals.
Using a result of Jarn\'\i k \cite{Jar36},
Laurent \cite{Lau09b} established the following statement. 

\proclaim Theorem 6.3.
Let $d$ and $n$ be integers with $n \ge 2$ and $0 \le d \le n-1$.
For every $\omega$ in $[(d+1)/(n-d), + \infty]$, there exists $\uxi$ such that
$\omega_{n, d} (\uxi) = \omega$. 
%that the spectrum of 
%$\om_{n,d}$ over ${\bf R}^n$ is equal to the whole
%interval $[(d+1)/(n-d), + \infty]$ and that 
Furthermore, $\om_{n,d} (\uxi) = (d+1)/(n-d)$ for almost all $\uxi$ in
${\bf R}^n$.

By means of the numbers $\xi_w$ 
defined in the proof of Theorem 5.1, we get some
information on the spectra of the exponents $w_{n, d}$.

\proclaim Theorem 6.4.
For $n \ge 2$ and $0 \le d \le n-1$, the spectrum of
$w_{n, d}$ contains the whole interval $[(n+d)/(n-d), + \infty]$ and
$w_{n, d} (\xi) = (d+1)/(n-d)$ for almost all
real numbers $\xi$.

Theorem 6.4 plainly includes the last assertion of Theorem 5.1 and 
the first assertion of Theorem 5.6. 

\pro
We follow the proof of the Corollary from \cite{Lau09b}, where it is
established that, for any $\lambda$ with $1/n \le \lambda \le + \infty$ 
and for any point $\utheta$ in ${\bf R}^n$ such that
$\om_{n, 0} (\utheta) = \lambda$ and $\om_{n, n-1} (\utheta) = n \lambda + n - 1$,
we have
$$
\om_{n, d} (\utheta) = {n \lambda + d \over n- d},
\quad (d=0, 1, \ldots , n-1).   \eqno (6.1)
$$
For $w\ge 2n-1$, the numbers $\xi_w$ defined in the
proof of Theorem 5.1 satisfy 
$$
n \lambda_n (\xi_w) =
w_n (\xi_w) - n + 1 = w - n +1,
$$
that is, 
$$ 
\om_{n, n-1} (\xi_w, \ldots , \xi_w^n)
= n \om_{n, 0} (\xi_w, \ldots , \xi_w^n) + n - 1.
$$
We then get from (6.1) that 
$$
w_{n, d} (\xi_w) = {n \lambda_n (\xi_w) + d \over n- d},
\quad (d=0, 1, \ldots , n-1).
$$
The first assertion of the theorem follows since 
$\lambda_n (\xi_w)$ takes every value between $1$ and $+ \infty$ 
as $w$ varies from $2n - 1$ to $+ \infty$.
The second assertion is an immediate consequence of (6.1)
and the fact that $n \lambda_n (\xi) =
w_n (\xi) - n + 1 = 1$ holds for almost every real number $\xi$.
\cqfd

\goodbreak

We conclude this section by mentioning that Laurent \cite{Lau09a}
determined the set of values taken by the quadruple of
functions $(\omega_{2,0},\omega_{2,1},\hom_{2,0},\hom_{2,1})$ at real points.

\vskip 5mm

\centerline{\bf 7. Parametric geometry of numbers}

\vskip 5mm

In 2009,  Schmidt and Summerer \cite{ScSu09,ScSu13a} introduced
a new theory, called the parametric geometry of numbers.
They studied the joint behaviour of the $n+1$ successive minima of 
certain one parameter families of convex bodies in $\R^{n+1}$, as a function of
the parameter. They further showed how their results allow them to recover 
classical inequalities relating various exponents of Diophantine approximation
attached to points in $\R^n$ and to find new ones. 
Some aspects of their theory have been simplified and completed by Roy \cite{Roy15a}, who
was then able to derive several spectacular results \cite{Roy15b,Roy15c}.

Let $n$ be a positive integer and $\uu$ be a non-zero vector in $\R^{n+1}$.
For each real number $Q\ge 1$, we form the convex body
$$
\cC_\uu(Q)
 = \{\ux\in\bR^{n+1} \,;\, \|\ux\|_2 \le 1, \ |\ux\cdot\uu| \le Q^{-1}\},
$$
where $| \cdot |$ denotes the scalar product and $\| \cdot \|_2$ the Euclidean norm. 
For $j=1,\dots,n+1$, we denote by
$\lambda_j\big(\cC_\uu(Q)\big)$
the $j$-th minimum of $\cC_\uu(Q)$, namely the smallest real number
$\lambda>0$ such that $\lambda \cC_\uu(Q)$ contains
at least $j$ linearly independent points of $\Z^{n+1}$.
Schmidt and Summerer \cite{ScSu13a} defined
$$
 L_{\uu,j}(q) = \log \lambda_j(\cC_\uu({\rm e}^q)),
 \quad
 q\ge 0,\ 1\le j\le n+1,
$$
and considered the map $\uL_\uu\colon[0,\infty)\to \bR^{n+1}$
given by
$$
 \uL_\uu(q) = (L_{\uu,1}(q),\dots,L_{\uu,n+1}(q)), 
 \quad  q\ge 0.
$$
They established many properties of this map. For instance, each of its
components $L_{\uu,j} : [0, + \infty) \to \R$ is continuous and piecewise
linear with slopes $0$ and $1$. Schmidt and Summerer showed that each 
function $\uL_\uu$ can be approximated up to bounded difference by
functions from a certain class, and Roy \cite{Roy15a} showed that the
same property holds within a simpler class.

For $j=1,\dots,n+1$, Roy \cite{Roy15b} also introduced 
$$
 \psibot_{j}(\uu)
  = \liminf_{q\to + \infty} {L_{\uu,1}(q)+\cdots+L_{\uu,j}(q) \over q}
\quad \hbox{and} \quad
 \psitop_{j}(\uu)
 = \limsup_{q\to + \infty} {L_{\uu,1}(q)+\cdots+L_{\uu,j}(q) \over q}.
$$
The following result, established in \cite{Roy15b}, 
connects these quantities to those from the previous section.

\proclaim Proposition 7.1.
Let $n$ be a positive integer and $\uu = (1, \xi_1, \ldots , \xi_n)$ be a vector in $\R^{n+1}$.  
For $j=0,\dots, n-1$, we have
$$
\omega_{n,j} ( (\xi_1, \ldots , \xi_n) )= {1 \over \psibot_{n-j}(\uu)}-1
\quad \hbox{and} \quad
 \hat{\omega}_{n,j} ( (\xi_1, \ldots , \xi_n) )= {1 \over \psitop_{n-j}(\uu)}-1.
$$

We quote below the main result of \cite{Roy15b} and observe
that it implies the first statement of Theorem 6.3.

\proclaim Theorem 7.2.
Let $n$ be a positive integer. 
For any $\omega_0,\ldots,\omega_{n-1}
\in [0, + \infty]$ satisfying $\omega_0\ge 1/n$ and
$$
{j\omega_{n,j} \over \omega_{n,j} +j+1}  \le \omega_{n, j-1}
\le {(n-j)\omega_{n,j} -1 \over n-j+1}, \quad
1\le j\le n-1,
$$
there exists a point $\uu\in\R^{n}$, whose coordinates are, together with $1$, linearly
independent over $\Q$, such that
$$ 
\omega_{n,j}(\uu)=\omega_j \quad \hbox{and} \quad  
\homega_{n,j} (\uu)={j+1 \over n-j}, \quad 0 \le j\le n-1.
$$

Furthermore, the point of view of parametric geometry of numbers has led 
Schmidt and Summerer to introduce the following exponents
of approximations.

\proclaim Definition 7.3.
Let $\utheta = (\xi_1, \ldots , \xi_n)$ be in $\R^n$. 
For $i=1, \ldots , n+1$, we denote by $\lambda_{i, n} (\utheta)$ (resp. $\hlambda_{i, n} (\utheta)$)
the supremum of the real numbers $\lambda$ such that the system of inequalities
$$
0 < |x_0| \le X, \quad \max_{1 \le j \le n} |x_0 \xi_j - x_j| \le X^{- \lambda}
$$
has $i$ linearly independent solutions $(x_0, x_1, \ldots , x_n)$ in $\Z^n$ for
arbitrarily large $X$ (resp. for every sufficiently large $X$). 

Schmidt and Summerer  \cite{ScSu13a} 
observed that $\lambda_{1, n} (\utheta) = \lambda_n (\utheta)$,
$\hlambda_{1, n} (\utheta) = \hlambda_n (\utheta)$,
$\lambda_{n+1, n} (\utheta) = 1 / \hw_n (\utheta)$ and
$\hlambda_{n+1, n} (\utheta) = 1 / w_n (\utheta)$, by 
Mahler's theorem on polar reciprocal bodies \cite{Mah39}. 

These exponents have been studied by Schleischitz \cite{Schl13,Schl15b}.

\vskip 5mm

\goodbreak

\centerline{\bf 8. Real numbers which are badly approximable by algebraic numbers}

\vskip 5mm

This short section is devoted to Problems 24 and 25 of \cite{BuLiv}, which were
solved by Badziahin and Velani \cite{BaVe14}, 
for $n = 2$, and by Beresnevich \cite{Ber15}, for $n \ge 3$.

Recall that an irrational real number $\xi$ is called a {\it badly approximable number}
if there exists a positive real number $c$ such that 
$$
|q x - p| > {c \over |q|}, \quad
\hbox{for every $p, q$ in $\Z$ with $q \not= 0$.}
$$
This notion can be extended as follows.

\proclaim Definition 8.1.
Let $n$ be a positive integer. A real number $\xi$ 
is called $n$-badly approximable if there exists a positive 
constant $c (\xi, n)$ such that
$$
|P(\xi)| \ge {c (\xi, n) \, H(P)^{-n}} 
\quad \hbox{for any integer polynomial $P(X)$ of degree $\le n$.}  
$$

Observe that it follows from Liouville's Inequality (Theorem 2.9) that, for any positive
integer $n$, any real algebraic number of degree $n+1$ is $n$-badly approximable. 

Davenport \cite{Dav64} asked whether there exist 
$2$-badly approximable transcendental real numbers. His question remained
open for nearly fifty years, until it was finally solved
by Badziahin and Velani \cite{BaVe14}, using an intricate construction
inspired by their proof \cite{BaPoVe11} of a conjecture of Schmidt.
Their result was subsequently extended a few years later by 
Beresnevich \cite{Ber15}, who confirmed the existence of $n$-badly approximable
transcendental real numbers, for every given positive integer $n$. For $n=2$, his proof differs
greatly from that of \cite{BaVe14}. 

A closely related problem deals with transcendental numbers 
badly approximable %\index{badly approximable number} 
by algebraic numbers of degree at most equal to some integer $n$. 
As well, it has been solved by Badziahin and Velani \cite{BaVe14},
for $n = 2$, and by Beresnevich \cite{Ber15}, for $n \ge 3$.

\proclaim Theorem 8.2. 
Let $n$ be a positive integer. 
There exist  transcendental real numbers $\xi$ which are
$n$-badly approximable and for which there exist positive 
real numbers $c_1(\xi, n)$ and $c_2(\xi, n)$ such that
$$
| \xi - \alpha | \ge {c_1(\xi, n) \, H(\alpha)^{-n-1}}, 
\quad \hbox{for any real algebraic number $\alpha$ of degree $\le n$,}  
$$
and
$$
| \xi - \alpha | \le {c_2(\xi, n) \, H(\alpha)^{-n-1}}, \quad
\hbox{for infinitely many real algebraic numbers
$\alpha$ of degree $\le n$.}  
$$
Moreover, the set of real numbers with this property
has full Hausdorff dimension.

Actually, Beresnevich \cite{Ber15} established that any intersection of finitely many of the sets 
arising in Theorem 8.2 has full Hausdorff dimension.

\vskip 5mm

\centerline{\bf 9. Open problems}

\vskip 5mm

In Chapter 10 of \cite{BuLiv} we listed several open questions.
As we have already seen above, some of them have now been solved. 
We gather below some of the still open problems mentioned in \cite{BuLiv}, and 
add a few supplementary ones. 

We begin with the conjecture of Wirsing \cite{Wir} dealing with
the approximation of real transcendental numbers by 
real algebraic numbers of bounded degree. This celebrated open problem
has motivated many works in this area. 

\proclaim Problem 9.1 (Wirsing's Conjecture).  
For any integer $n\ge 1$ and for any
real transcendental number $\xi$, we have $w_n^*(\xi) \ge n$.

We may even ask for a stronger version of Wirsing's conjecture, namely whether,
for any positive integer $n$ and any
real transcendental number $\xi$, there exist 
a constant $c(\xi, n)$ and
infinitely many real algebraic numbers $\alpha$ of degree less
than or equal to $n$ such that
$$
|\xi - \alpha| \le c(\xi, n) \, H(\alpha)^{-n-1}.
$$
Davenport and Schmidt \cite{DaSc67} gave a positive answer
to this question in the case $n=2$ (and Theorem 8.2 implies that their result is best possible up to
the value of the numerical constant), 
but we do not know if we can fix the exact degree of the approximants instead of just 
an upper bound for it. 

\proclaim Problem 9.2. 
For any integer $n\ge 2$ and any 
real transcendental number $\xi$, there exist 
a constant $c(\xi, n)$ and
infinitely many real algebraic numbers $\alpha$ of degree exactly $n$ such that
$$
|\xi - \alpha| \le c(\xi, n) \, H(\alpha)^{-n-1}.
$$

Results of Roy \cite{Roy03b,Roy04a} could
speak in favour of the existence of transcendental numbers which do not
satisfy the conclusion of Problem 9.2, even for $n=2$.

The next problem was called the `main problem' in \cite{BuLiv}.

\proclaim Problem 9.3. 
Let $(w_n)_{n \ge 1}$ and $(w_n^*)_{n \ge 1}$
be two non-decreasing sequences in $[1, + \infty]$ such that
$$
n \le w_n^* \le w_n \le w_n^* + n - 1, \quad {\hbox{for any $n \ge 1$}}.
$$
Then there exists a real transcendental number $\xi$ such that
$$
w_n(\xi) = w_n \quad {\it and} \quad w_n^*(\xi) = w_n^*
\quad {\hbox{for any $n \ge 1$}}.
$$

Since Problem 9.3 does not take the exponents $\lambda_n$ into account, we propose a
more general formulation.

\proclaim Problem 9.4.
For $n \ge 2$, determine the joint spectrum of the triple of exponents
$(w_n, w_n^*, \lambda_n)$, that is, the set of triples
$(w_n(\xi), w_n^*(\xi), \lambda_n(\xi))$, when $\xi$ goes through the set of 
transcendental real numbers.

Theorem 3.2 shows how the exponents $\lambda_n$ are 
related to the exponents $w_n$ by means of a transference theorem.

\proclaim Problem 9.5. 
Let $n \ge 2$ be an integer, $\lambda_n \ge 1/ n$
and $w_n \ge n$ be real numbers satisfying
$$
{w_n \over (n-1) w_n + n} \le 
\lambda_n \le {w_n -n + 1 \over  n}.
$$
There exist real numbers $\xi$ such that $w_n(\xi) = w_n$ and
$\lambda_n (\xi) = \lambda_n$.

Corollary 2.7 asserts that $w_n^*(\xi)=n$ holds if $w_n(\xi) = n$,
but the converse is an open question.

\proclaim Problem 9.6. 
For any positive integer $n$, we have $w_n(\xi) = n$
if $w_n^*(\xi) = n$.

We now turn to uniform exponents. 

\proclaim Problem 9.7. 
For any integer $n\ge 3$ and any
real transcendental number $\xi$, we have $\hw_n(\xi) =  n$.
At least, obtain a better upper bound than $\hw_n (\xi) \le 2n - 1$.

Approximation by algebraic integers (that is, by algebraic numbers, whose
minimal defining polynomial over $\Z$ is monic) has been first studied 
by Davenport and Schmidt \cite{DaSc69}. Roy \cite{Roy03b} 
proved that there exist real numbers $\xi$ which are very badly approximable by
quadratic integers, in the sense that there exists $c > 0$ such that
$$
|\xi - \alpha| > c H(\alpha)^{-(1 + \sqrt{5})/2},
$$
for every real quadratic number $\alpha$. 

The next problem is the analogue of Wirsing's conjecture for 
approximation by algebraic integers.

\proclaim Problem 9.8.
For any integer $n\ge 4$, any positive real number $\eps$, and any
real transcendental number $\xi$, there exist 
a constant $c(\xi, n, \eps)$ and
infinitely many real algebraic integers $\alpha$ of degree less
than or equal to $n$ such that
$$
|\xi - \alpha| \le c(\xi, n, \eps) \, H(\alpha)^{-n + \eps}.
$$

In view of auxiliary results from \cite{DaSc69}, the answer to Problem 9.8
is positive for some integer $n \ge 4$ if one can prove that any
real transcendental number $\xi$ satisfies $\hw_{n-1} (\xi) = n-1$. 

Despite the recent progress made in \cite{BaVe14,Ber15}, the 
following problem remains open.

\proclaim Problem 9.9.  
There exist a real transcendental number $\xi$ and a sequence $(c(\xi, n))_{n \ge 1}$
of positive real numbers such that
$$
|P(\xi)| \ge {c(\xi, n) \, H(P)^{-n}} 
\quad \hbox{for any integer $n$ and any polynomial $P(X)$ of degree $\le n$.}   
$$

It is likely that the answer to Problem 9.9 is positive and that, moreover,
the set of real numbers $\xi$ with this property has full
Hausdorff dimension.

We continue with a problem proposed by Schleischitz \cite{Schl15a}, 
which corrects and refines
a problem posed in \cite{Bu10a}. 

\proclaim Problem 9.10.
Let $m, n$ be integers with $1 \le n \le m$. Does the inequality
$$
\lambda_m (\xi) \ge {n \lambda_n (\xi) - m + n \over m}
$$
hold for every transcendental real number $\xi$?

The next problem extends a question 
posed by Beresnevich, Dickinson, and Velani \cite{BDV} in the case
of (simultaneous) rational approximation.

\proclaim Problem 9.11. 
Let $n$ be a positive integer
and let $\tau > 1$ be real. Is the set of real numbers $\xi$ for which there
exists a positive constant $c (\xi)$ such that
$$
|\xi - \alpha| \le H(\alpha)^{-\tau (n+1)},  \quad
\hbox{for infinitely many algebraic numbers $\alpha$ of degree $\le n$,}
$$
and
$$
|\xi - \alpha| \ge c (\xi) \, H(\alpha)^{-\tau (n+1)}, \quad
\hbox{for every algebraic number $\alpha$ of degree $\le n$,}
$$
non-empty? If yes, determine its Hausdorff dimension.

Problem 9.11 has been solved \cite{Bu03} when $n=1$.
One may also replace the approximation functions $x \mapsto x^{-\tau (n+1)}$ 
by more general non-increasing functions $\Psi$.

Problems 9.12 and 9.13 deal with metrical results. 

\proclaim Problem 9.12. 
Let $n \ge 2$ be an integer. Let $\lambda_n$ be
a real number with $\lambda_n \ge 1/n$. Determine the Hausdorff dimension
of the set of real numbers $\xi$ such that $\lambda_n (\xi) = \lambda_n$.

\proclaim Problem 9.13. 
Determine the Hausdorff dimension
of the set of real numbers $\xi$ such that $\hw_2 (\xi) > 2$
(\resp $\hw^*_2 (\xi) > 2$).

We continue with two specific questions on exponents of
uniform approximation. 

\proclaim Problem 9.14.
Does there exist $\xi$ such that $\hw_2 (\xi) > 2$ and
$\hw_2 (\xi) > \hw_2^* (\xi)$?

\proclaim Problem 9.15.
The spectrum of $\hw_2^*$ includes the interval $[1, 2]$. 

Let $m \ge 2$ be an integer.
According to LeVeque \cite{LeV53}, a real number $\xi$ is a $U_m$-number if
$w_m (\xi)$ is infinite and $w_{m-1} (\xi)$ is finite. Furthermore, the $U_1$-numbers are
precisely the Liouville numbers. 

It is proved in \cite{AlAvBu03} (see also Section 7.6 of \cite{BuLiv}) that, for any
integer $m \ge 2$, there exist uncountably many real $U_m$-numbers $\xi$ 
with  
$$
w_n^* (\xi) \le m + n - 1, \quad \hbox{for $n = 1, \ldots , m-1$}.    \eqno (9.1)
$$
Schmidt \cite{Schm06} showed that $w_n^* (\xi)$ can be replaced by
$w_n (\xi)$ in (9.1).
%, thereby he obtained a stronger result.

\proclaim Problem 9.16.
Let $m \ge 2$ be an integer. There exist real $U_m$-numbers $\xi$ 
satisfying $w_n (\xi) = n$, for $n = 1, \ldots , m-1$. 

Corollary 5.4 shows that the values taken by our exponents
of approximation at $U_1$-numbers are known. 

\proclaim Problem 9.17. 
Let $m$ and $n$ be integers with $n \ge m \ge 2$. 
Study the values taken by the exponents of approximation $\lambda_n$,
$\hw_n^*, \ldots$ at $U_m$-numbers.

Among many questions concerning the exponents $w_{n,d}$ and $\hw_{n,d}$
defined in Section 6, let us point out the following three ones.

\proclaim Problem 9.18.
Let $d$ and $n$ be integers with $n \ge 2$ and $0 \le d \le n-1$. 
Find a real number $C_{n,d}$ as small as possible such that every
transcendental real number $\xi$ satisfies $\hw_{n,d} (\xi) \le C_{n,d}$.

\proclaim Problem 9.19. 
Let $d$ and $n$ be integers with $n \ge 2$ and $0 \le d \le n-1$. 
Is the spectrum of the function $w_{n, d}$
equal to $[(d+1)/(n-d), + \infty]$ ?

A positive answer of Problem 9.19 would (probably) follow from the resolution of
the next problem. 

\proclaim Problem 9.20.
Let $d$ and $n$ be integers with $n \ge 2$ and $0 \le d \le n-1$. 
Let $w$ be a real number satisfying $w > (d+1) / (n-d)$. 
Determine the Hausdorff dimension of the sets
$$
\{\uxi \in \R^n : \omega_{n,d} (\uxi) \ge w\}
$$
and
$$
\{\xi \in \R : w_{n,d} (\xi) \ge w\}.
$$

Throughout this survey, we have focus our attention on approximation to
real numbers. However, we may as well
consider approximation to complex numbers or to $p$-adic
numbers; see the references given in Chapter 9 of \cite{BuLiv} and the  
works \cite{BuEv09,Zelo,Bu10b,BuBuDiO,Bel13,Ke13,BuPe15}.

There are as well several recent papers on uniform Diophantine approximation on
curves in $\R^2$; see \cite{LoRoy12,Roy13,Bat15}.

\vskip 8mm

%\vfill\eject

\centerline{\bf References}

\vskip 5mm

\beginthebibliography{999}

\bibitem{AlAvBu03}
K. Alnia\c cik, Y. Avci, and Y. Bugeaud,
{\it On $U_m$ numbers with small transendence measure},
Acta Math. Hungar. 99 (2003), 271--277.

\bibitem{ArRo}
B. Arbour and D. Roy,
{\it A Gel'fond type criterion in degree two},
Acta Arith. 111 (2004), 97--103.

\bibitem{BaPoVe11}
D. Badziahin, A. Pollington, and S. Velani,
{\it On a problem in simultaneous diophantine
approximation: Schmidt's conjecture},
Ann. of Math. 174 (2011), 1837--1883.

\bibitem{BaVe14}
D. Badziahin and S. Velani,
{\it Badly approximable points on planar curves and a problem of Davenport},
Math. Ann.  359 (2014), 969--1023.

\bibitem{BaSc70}
A. Baker and W. M. Schmidt,
{\it Diophantine approximation and Hausdorff dimension},
Proc. London Math. Soc. 21 (1970), 1--11.

\bibitem{Bat15}
G. Batzaya,
{\it On simultaneous approximation to powers of a real number by rational numbers},
J. Number Theory 147 (2015), 141--155.

\bibitem{Bel13}
P. Bel,
{\it Approximation simultan\'ee d'un nombre $v$-adique et de son carr\'e 
par des nombres alg\'ebriques}, 
J. Number Theory 133 (2013),  3362--3380.

\bibitem{Ber15}
V. Beresnevich,  
{\it Badly approximable points on manifolds}.
Preprint.

\bibitem{BDV}
V. V. Beresnevich, H. Dickinson, and S. L. Velani,
{\it Sets of exact `logarithmic order' in the theory of Diophantine approximation}, 
Math. Ann. 321 (2001), 253--273.

\bibitem{BeDiVe07}
V. Beresnevich, D. Dickinson and S. L. Velani,
{\it Diophantine approximation on planar curves and the distribution of rational points}, 
with an appendix by R.C. Vaughan: `Sums of two squares near perfect squares',
Ann. of Math. 166 (2007), 367--426.

\bibitem{Ber83}
V. I. Bernik,
{\it Application of the Hausdorff dimension in the theory of Diophantine approximations}, 
Acta Arith. 42 (1983), 219--253 (in Russian).
English transl. in Amer. Math. Soc. Transl. 140 (1988), 15--44.

\bibitem{BeTi}
V. I. Bernik and K. Tishchenko,
{\it Integral polynomials with an overfall of the coefficient values and Wirsing's problem}, 
Dokl. Akad. Nauk Belarusi 37 (1993), no. 5, 9--11 (in Russian).

\bibitem{BuBuDiO}
N. Budarina, Y. Bugeaud, D. Dickinson, and H. O'Donnell,  
{\it On simultaneous rational approximation to a $p$-adic number and its integral powers},
Proc. Edinb. Math. Soc. 54 (2011), 599--612.

\bibitem{Bu03}
Y. Bugeaud,
{\it Sets of exact approximation order by rational numbers},
Math. Ann. 327 (2003), 171--190.

\bibitem{BuLiv}
Y. Bugeaud,
Approximation by algebraic numbers,
Cambridge Tracts in Mathematics, Cambridge, 2004.

\bibitem{Bu10a}
Y. Bugeaud,
{\it On simultaneous rational approximation
to a real number and its integral powers},
Ann. Inst. Fourier (Grenoble) 60 (2010), 2165--2182.

\bibitem{Bu10b}
Y. Bugeaud,
{\it On simultaneous uniform approximation to a $p$-adic number and its square},
Proc. Amer. Math. Soc. 138 (2010), 3821--3826.

\bibitem{Bu12b}
Y. Bugeaud,
{\it Continued fractions with low complexity: 
Transcendence measures and quadratic approximation},
Compos. Math. 148 (2012), 718--750.

\bibitem{BuDu11}
Y. Bugeaud and A. Dujella,
{\it Root separation for irreducible integer polynomials},
Bull. Lond. Math. Soc. 43 (2011), 1239--1244.

\bibitem{BuEv09}
Y. Bugeaud and J.-H. Evertse,
{\it Approximation of complex algebraic numbers 
by algebraic numbers of bounded degree},
Ann. Scuola Normale Superiore di Pisa 8 (2009), 333--368.

\bibitem{BuLa05a}
Y. Bugeaud and M. Laurent,
{\it Exponents of Diophantine Approximation and Sturmian Continued Fractions},
Ann. Inst. Fourier (Grenoble) 55 (2005), 773--804.

\bibitem{BuLa05b}
Y. Bugeaud and M. Laurent,
{\it Exponents of homogeneous and inhomogeneous Diophantine Approximation},
Moscow Math. J. 5 (2005), 747--766.

\bibitem{BuLa07}
Y. Bugeaud and M. Laurent,
{\it Exponents of Diophantine approximation}.
In: Diophantine Geometry Proceedings,
Scuola Normale Superiore Pisa, Ser. CRM, vol. 4, 2007, 101--121.

\bibitem{BuLa10}
Y. Bugeaud and M. Laurent,
{\it On transfer inequalities in Diophantine approximation, II},
Math. Z.  265  (2010), 249--262. 

\bibitem{BuPe15}
Y. Bugeaud and T. Pejkovic,
{\it Quadratic approximation in $\Q_p$}, 
Intern. J. Number Theory 11 (2015), 193--209.

\bibitem{BuTeu}
Y. Bugeaud et O. Teuli\'e,
{\it Approximation d'un nombre r\'eel par des nombres alg\'ebriques de degr\'e donn\'e}, 
Acta Arith. 93 (2000), 77--86.

\bibitem{Cas}
J. W. S. Cassels,
An introduction to Diophantine Approximation,
Cambridge Tracts in Math. and Math. Phys., vol. 99, Cambridge
University Press, 1957.

\bibitem{Dav64}
H. Davenport,
{\it A note on Diophantine approximation. II}, 
Mathematika 11 (1964), 50--58.

\bibitem{DaSc67}
H. Davenport and W. M. Schmidt,
{\it Approximation to real numbers by quadratic irrationals}, 
Acta Arith. 13 (1967), 169--176.

\bibitem{DaSc68} 
H. Davenport and W. M. Schmidt,
{\it A theorem on linear forms}, 
Acta Arith. 14 (1967/1968), 209--223.

\bibitem{DaSc69}
{H. Davenport and W. M. Schmidt},
{\it Approximation to real numbers by algebraic integers}, 
Acta Arith. {15} (1969), 393--416.

\bibitem{Fisc04}
S. Fischler,
{\it Spectres pour l'approximation d'un nombre r\'eel et de son carr\'e},
C. R. Acad. Sci. Paris 339 (2004), 679--682.

\bibitem{Fisc06}
S. Fischler,
{\it Palindromic prefixes and episturmian words},
J. Combin. Theory Ser. A  113 (2006), 1281--1304.

\bibitem{Fisc07}
S. Fischler,
{\it Palindromic prefixes and Diophantine approximation},
Monatsh. Math.  151 (2007),  11--37.

\bibitem{Ger12a}
O. N. German,
{\it On Diophantine exponents and Khintchine's transference principle}, 
Mosc. J. Comb. Number Theory 2 (2012), 22--51.

\bibitem{Jar36}
V. Jarn\'\i k,
{\it \"Uber einen Satz von A. Khintchine, 2. Mitteilung},
Acta Arith. 2 (1936), 1--22.

\bibitem{Jar38}
V. Jarn\'\i k,
{\it Zum Khintchineschen ``\"Ubertragungssatz''},
Trav. Inst. Math. Tbilissi 3 (1938), 193--212.

\bibitem{Jar50}
V. Jarn\'\i k,
{\it Une remarque sur les approximations diophantiennes lin\'eaires},
Acta Sci. Math. Szeged 12 (1950), 82--86.

\bibitem{Jar54}
V. Jarn\'\i k,
{\it Contribution \`a la th\'eorie des approximations diophantiennes
lin\'eaires et homog\`enes},
Czechoslovak Math. J. 4 (1954), 330--353 (in Russian, French summary).

\bibitem{Ke13}
G. Kekec,
{\it On Mahler's $p$-adic $U_m$-numbers}, 
Bull. Aust. Math. Soc. 88 (2013), 44--50.

\bibitem{Kh26b}
A. Ya. Khintchine,
{\it \"Uber eine Klasse linearer diophantischer Approximationen},
Rendiconti Circ. Mat. Palermo 50 (1926), 170--195.

\bibitem{Kh48b}
A. Ya. Khintchine,
{\it On some applications of the method of the additional variable}, 
Uspehi Matem. Nauk 3, (1948), 188--200 (in Russian). English translation:
Amer. Math. Soc. Translation (1950), no. 18, 14 pp.

\bibitem{Ko39}
J. F. Koksma,
{\it \"Uber die Mahlersche Klasseneinteilung der transzendenten Zahlen
und die Approximation komplexer Zahlen durch algebraische Zahlen},
Monatsh. Math. Phys. 48 (1939), 176--189.

\bibitem{Lau03}
M. Laurent,
{\it Simultaneous rational approximation to the successive powers of a real number},
Indag. Math. 11 (2003), 45--53.

\bibitem{Lau09a}
M. Laurent,
{\it Exponents of Diophantine approximation in dimension two},
Canad. J. Math.  61  (2009), 165--189. 

\bibitem{Lau09b}
M. Laurent,
{\it On transfer inequalities in Diophantine approximation}.  
In: Analytic number theory,  306--314, Cambridge Univ. Press, Cambridge, 2009.

\bibitem{LeV53}
W. J. LeVeque,
{\it On Mahler's $U$-numbers},
J. London Math. Soc. 28 (1953), 220--229.

\bibitem{LoRoy12}
S.  Lozier and D. Roy, 
{\it Simultaneous approximation to a real number and to its cube by rational numbers}, 
Acta Arith. 156 (2012), 39--73.

\bibitem{Mah32}
K. Mahler,
{\it Zur Approximation der Exponentialfunktionen und des
Logarithmus. I, II},
J. reine angew. Math. 166 (1932), 118--150.

\bibitem{Mah39}
K. Mahler,
{\it Ein \"Ubertreibungsprinzip f\"ur konvexe K\"orper},
\v Casopis P\v est. Mat. Fyz. 68 (1939), 93--102. 

\bibitem{Mar15}
A. Marnat,
{\it About Jarn\'\i k-type relations in higher dimension.}
In preparation.

\bibitem{Mosh10}
N. G. Moshchevitin, 
{\it Singular Diophantine systems of A. Ya. Khinchin and their application},
Uspekhi Mat. Nauk 65 (2010), 43--126; English translation in
Russian Math. Surveys 65 (2010), 433--511.

\bibitem{Mosh12}
N. Moshchevitin, 
{\it Exponents for three-dimensional simultaneous Diophantine approximations}, 
Czechoslovak Math. J. 62 (2012), 127--137.

%\bibitem{Mosh13}
%N. Moshchevitin, 
%{\it Diophantine exponents for systems of linear forms in two variables},
%Acta Sci. Math. (Szeged)  79  (2013),  347--367.

\bibitem{Mosh14}
N. Moshchevitin, 
{\it A note on two linear forms},
Acta Arith. 162 (2014), 43--50.

\bibitem{Roy03a}
D. Roy,
{\it Approximation simultan\'ee d'un nombre et son carr\'e},
C. R. Acad. Sci. Paris 336 (2003), 1--6.  

\bibitem{Roy03b}
D. Roy,
{\it Approximation to real numbers by cubic algebraic numbers, II},
Ann. of Math. 158 (2003),  1081--1087.

\bibitem{Roy04a}
D. Roy,
{\it Approximation to real numbers by cubic algebraic numbers, I},
Proc. London Math. Soc. 88 (2004), 42--62.

\bibitem{Roy04b}
D. Roy,
{\it Diophantine approximation in small degree},
Number Theory, 269--285, CRM Proc. Lecture Notes 36,
Amer. Math. Soc., Providence, RI, 2004.

\bibitem{Roy07}
D. Roy,
{\it On two exponents of approximation related to a real number and
its square},
Canad. J. Math. 59 (2007), 211--224. 

\bibitem{Roy08a}
D. Roy,
{\it On the continued fraction expansion of a class of numbers}.
In: Diophantine approximation, Festschrift for Wolfgang Schmidt, 
Developments in Math. vol. 16, Eds: H. P. Schlickewei, K. Schmidt and R. Tichy, 
Springer-Verlag, 2008, 347--361.

\bibitem{Roy08b}
D. Roy,
{\it On simultaneous rational approximations to a real number, its square, and its cube}, 
Acta Arith. 133 (2008), 185--197. 

\bibitem{Roy13}
D. Roy,
{\it Rational approximation to real points on conics},
Ann. Inst. Fourier (Grenoble) 63 (2013), 2331--2348. 

\bibitem{Roy15a}
D. Roy,
{\it On Schmidt and Summerer parametric geometry of numbers},
Ann. of Math. 
To appear.

\bibitem{Roy15b}
D. Roy,
{\it Spectrum of the exponents of best rational approximation}.
Preprint.

\bibitem{Roy15c}
D. Roy,
{\it Construction of points realizing the regular systems 
of Wolfgang Schmidt and Leonard Summerer},
J. Th\'eor. Nombres Bordeaux.
To appear.

\bibitem{Schl13}
J. Schleischitz,
{\it Diophantine approximation and special Liouville numbers},
 Commun. Math. 21 (2013), 39--76.

\bibitem{Schl14}
J. Schleischitz,
{\it Two estimates concerning classical Diophantine approximation constants},
Publ. Math. Debrecen  84 (2014), 415--437.

\bibitem{Schl15a}
J. Schleischitz,
{\it On the spectrum of Diophantine approximation constants}.
Preprint.  arXiv:1409.1472.

\bibitem{Schl15b}
J. Schleischitz,
{\it On approximation constants for Liouville numbers}.
Preprint.  arXiv:1409.1396.

\bibitem{Schm67}
W. M. Schmidt,
{ \it On heights of algebraic subspaces and diophantine approximations},
Ann. of Math. 85 (1967), 430--472.

\bibitem{SchmLN}
W. M. Schmidt,
Diophantine Approximation. 
Lecture Notes in Math. {785}, Springer, Berlin, 1980.

\bibitem{Schm06}
W. M. Schmidt,
{\it Mahler and Koksma classification of points in $\R^n$ and $\C^n$}, 
Funct. Approx. Comment. Math. 35 (2006), 307--319.

\bibitem{ScSu09}
W. M. Schmidt and L. Summerer,
{\it Parametric geometry of numbers and applications}, 
Acta Arith.  140  (2009),  67--91.

\bibitem{ScSu13a}
W. M. Schmidt and L. Summerer,
{\it Diophantine approximation and parametric geometry of numbers}, 
Monatsh. Math. 169 (2013), 51--104.

\bibitem{ScSu13b}
W. M. Schmidt and L. Summerer,
{\it Simultaneous approximation to three numbers},
Moscow J. Comb. Number Th. 3 (2013), 84--107.

\bibitem{SprLiv}
{V. G. Sprind\v zuk},
Mahler's problem in metric number theory. 
Izdat. ``Nauka i Tehnika'', Minsk, 1967 (in Russian).
English translation by B. Volkmann, Translations of Mathematical Monographs,
Vol. 25, American Mathematical Society, Providence, R.I., 1969.

\bibitem{Tish00}
K. I. Tishchenko,
{\it On approximation to real numbers by algebraic numbers},
Acta Arith. 94 (2000), 1--24.

\bibitem{Tish07}
K. I. Tsishchanka, 
{\it On approximation of real numbers by algebraic numbers of bounded degree},
J. Number Theory 123 (2007), 290--314. 

\bibitem{VaVe06}
R. C. Vaughan and S. Velani,
{\it Diophantine approximation on planar curves: the convergence theory},
Invent. Math. 166 (2006), 103--124.

\bibitem{Wir}
{E. Wirsing},
{\it Approximation mit algebraischen Zahlen beschr\"ankten Grades}, 
J. reine angew. Math. {206} (1961), 67--77.

\bibitem{Zelo}
D. Zelo,
Simultaneous approximation to real and $p$-adic numbers.
PhD thesis. Univ. Ottawa. 2009. {\tt arXiv:0903.0086}.

\endthebibliography

\vskip1cm

\noindent Yann Bugeaud  

\noindent Universit\'e de Strasbourg

\noindent Math\'ematiques

\noindent 7, rue Ren\'e Descartes      

\noindent 67084 STRASBOURG  (FRANCE)

\vskip2mm

\noindent {\tt bugeaud@math.unistra.fr}

\bye